\documentclass[reqno,11pt]{amsart}
\usepackage{amsmath,graphicx,amssymb}
\usepackage[all]{xy}
\usepackage{epic}
\usepackage{eepic}
\newtheorem{proposition}{Proposition}%[section] % 2nd argument is what is printed
\newtheorem{cor}[proposition]{Corollary}
\newtheorem{lemma}[proposition]{Lemma}
\newtheorem{remark}[proposition]{Remark}

\numberwithin{equation}{section}

\renewcommand{\AA}{\mathfrak A}
\newcommand{\Z}{\mathbb Z}
\newcommand{\N}{\mathbb N}
\newcommand{\R}{\mathbb R}
\newcommand{\C}{\mathbb C}
\newcommand{\Q}{\mathbb Q}
\newcommand{\GG}{\mathcal G}
\newcommand{\A}{\mathfrak A}
\newcommand{\CCC}{\mathfrak C}
\newcommand{\VV}{\mathcal V}
\newcommand{\EE}{\mathcal E}
\newcommand{\G}{\mathbb G}
\newcommand{\F}{\mathbb F}
\renewcommand{\H}{\mathbb H}
\newcommand{\I}{\mathbb I}

\newcommand{\CC}{\mathfrak F}

\newcommand{\RR}{\mathcal R}
\newcommand{\BB}{\mathfrak B}
\newcommand{\PP}{\mathcal P}
\newcommand{\M}{\mathbb M}

\newcommand{\K}{\mathbb K}
\newcommand{\II}{\mathfrak I}

\newcommand{\OO}{\mathcal O}

\newcommand{\we}{\widetilde{e}}
\newcommand{\wf}{\widetilde{f}}
\newcommand{\Prim}{\operatorname{Prim}}

\newcommand{\TT}{\mathcal T}

\title[An AF algebra associated with the Farey tessellation]{An AF algebra associated with the \\ Farey tessellation}
\author{Florin P. Boca}

\address{Department of Mathematics, University of Illinois, 1409 W. Green Street, Urbana, IL 61801, USA}

\address{Institute of Mathematics ``Simion Stoilow" of the Romanian Academy,
P.O. Box 1-764, RO-014700 Bucharest, Romania}

\address{E-mail: fboca@math.uiuc.edu}


\date{June 20, 2008}

\thanks{}

\begin{document}

%\maketitle

\begin{abstract}
To the Farey tessellation of the upper half-plane we associate an
AF algebra $\AA$ encoding the cutting sequences that define
vertical geodesics. The Effros-Shen AF algebras arise as quotients
of $\AA$. Using the path algebra model for AF algebras we
construct, for each $\tau \in (0,\frac{1}{4}]$, projections
$(E_n)$ in $\AA$ such that $E_n E_{n\pm 1}E_n \leq \tau E_n$.
\end{abstract}

\maketitle

\tableofcontents

\section*{Introduction}
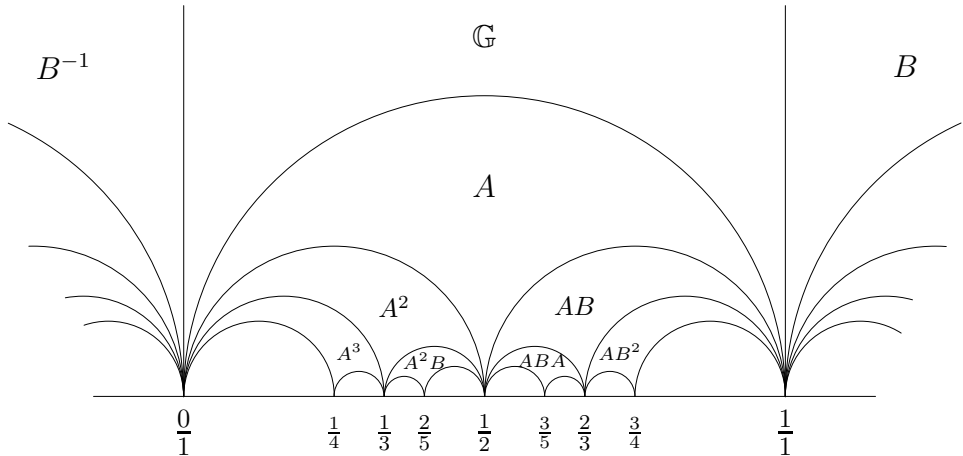
\begin{figure}[ht]
\begin{center}
\unitlength 0.4mm
\begin{picture}(0,140)(0,-12)

\path(-100,130)(-100,0) \path(100,130)(100,0) \path(-130,0)(130,0)
\put(0,120){\makebox(0,0){{\large $\G$}}}
\put(140,110){\makebox(0,0){{\large $B$}}}
\put(-140,110){\makebox(0,0){{\large $B^{-1}$}}}
\put(0,70){\makebox(0,0){{\large $A$}}}

\put(200,0){\arc{200}{3.14}{-2}} \put(0,0){\arc{200}{3.14}{0}}
\put(-200,0){\arc{200}{-1.14}{0}}

\put(50,0){\arc{100}{3.14}{0}} \put(-50,0){\arc{100}{3.14}{0}}
\put(150,0){\arc{100}{3.14}{-1.5}}
\put(-150,0){\arc{100}{-1.6}{0}}

\put(16.666,0){\arc{33.333}{3.14}{0}}
\put(66.666,0){\arc{66.666}{3.14}{0}}
\put(-16.666,0){\arc{33.333}{3.14}{0}}
\put(-66.666,0){\arc{66.666}{3.14}{0}}

\put(133.333,0){\arc{66.666}{3.14}{-1.3}}
\put(125,0){\arc{50}{3.14}{-1}}
\put(-133.333,0){\arc{66.666}{-1.75}{0}}
\put(-125,0){\arc{50}{-1.9}{0}}

\put(41.666,0){\arc{16.666}{3.14}{0}}
\put(75,0){\arc{50}{3.14}{0}} \put(10,0){\arc{20}{3.14}{0}}
\put(26.666,0){\arc{13.333}{3.14}{0}}

\put(-41.666,0){\arc{16.666}{3.14}{0}}
\put(-75,0){\arc{50}{3.14}{0}} \put(-10,0){\arc{20}{3.14}{0}}
\put(-26.666,0){\arc{13.333}{3.14}{0}}

\put(-100,-12){\makebox(0,0){{\Large $\frac{0}{1}$}}}
\put(100,-12){\makebox(0,0){{\Large $\frac{1}{1}$}}}
\put(0,-12){\makebox(0,0){{\large $\frac{1}{2}$}}}
\put(-33.33,-12){\makebox(0,0){{$\frac{1}{3}$}}}
\put(-30,30){\makebox(0,0){{\small $A^2$}}}
\put(30,30){\makebox(0,0){{\small $AB$}}}
\put(45,15){\makebox(0,0){{\tiny $AB^2$}}}
\put(-45,15){\makebox(0,0){{\tiny $A^3$}}}
\put(-20,12){\makebox(0,0){{\tiny $A^2B$}}}
\put(19,12){\makebox(0,0){{\tiny $ABA$}}}
\put(33.33,-12){\makebox(0,0){{$\frac{2}{3}$}}}
\put(-50,-12){\makebox(0,0){{\small $\frac{1}{4}$}}}
\put(50,-12){\makebox(0,0){{\small $\frac{3}{4}$}}}
\put(-20,-12){\makebox(0,0){{$\frac{2}{5}$}}}
\put(20,-12){\makebox(0,0){{\small $\frac{3}{5}$}}}

\end{picture}
\end{center}
\caption{The Farey tessellation} \label{Figure1}
\end{figure}

The semigroup ${\mathfrak S}$ generated by the matrices $A=\left[
\begin{smallmatrix} 1 & 0
\\ 1 & 1 \end{smallmatrix}\right]$ and $B=\left[ \begin{smallmatrix} 1 & 1
\\ 0 & 1 \end{smallmatrix}\right]$ is isomorphic to $\F_2^+$, the free
semigroup on two generators. This fact, intimately connected to
the continued fraction algorithm, can be visualized by means of
the \emph{Farey tessellation} $\{ g\G : g\in {\mathfrak S}\}$ of
$\H$ depicted in Figure \ref{Figure1}, where $\G=\big\{ 0\leq \Re
z\leq 1:\big| z-\frac{1}{2}\big|\geq \frac{1}{2}\big\}$ (cf.,
e.g., \cite{Se}).

The half-strip $0\leq \Re z\leq 1$, $\Im z >0$, is tessellated
precisely by the images of $\G$ under matrices from the set
\begin{equation*}
{\mathfrak S}_* =\{ I\} \cup \left\{ \left[ \begin{smallmatrix} a
& b \\ c & d \end{smallmatrix}\right]\in SL_2(\Z): 0\leq a \leq
c,\ 0\leq b\leq d \right\}.
\end{equation*}

By suspending the cusps in this tessellation (which correspond to
rational numbers in $[0,1]$) with appropriate (infinite)
multiplicities, one gets the diagram $\GG$ from Figure
\ref{Figure2} (cf. \cite{Kn2}). This diagram reflects both the
elementary mediant construction, that produces from a pair $(
\frac{p}{q},\frac{p'}{q'})$ of rational numbers with $p'q-pq'=1$
the new pairs $(\frac{p}{q},\frac{p+p'}{q+q'})$ and $(
\frac{p+p'}{q+q'},\frac{p'}{q'})$ with the same property, and the
geometry of the continued fraction algorithm. As in the case of
the Pascal triangle, in $\GG$ one writes the sum of the
denominators of two neighbors from the same floor into the next
floor of the diagram. One keeps, however, a copy of each
denominator at the next floor. For this reason, such a diagram was
called the \emph{Pascal triangle with memory} \cite{Kn1}. There is
a remarkable one-to-one correspondence between the integer
solutions of the equation $ad-bc=1$ with $0\leq a\leq c$, $0\leq
b\leq d$, and the rational labels of two neighbors at the same
floor in $\GG$, acquired by the mediant construction and by
keeping each label at the next floor in the diagram.

The thrust of this paper is the remark that, by regarding $\GG$ as
a Bratteli diagram, one gets an AF algebra $\AA=\varinjlim \AA_n$
with interesting properties. This algebra is closely related with
the \emph{Effros--Shen AF algebras} \cite{ES,PV} which we show to
arise as primitive quotients of $\AA$. The primitive ideal space
$\Prim \AA$ is identified with the disjoint union of the
irrational numbers in $[0,1]$ and three copies of the rational
ones, except for the endpoints $0$ and $1$ which are represented
by only two copies.

In \cite{Bra2} it was shown that any separable abelian
$C^*$-algebra $\mathfrak{Z}$ is the center $Z ({\mathcal A})$ of
an AF algebra ${\mathcal A}$. The AF algebra $\AA$ can actually be
retrieved from that abstract construction by embedding
$\mathfrak{Z}=C[0,1]$ into the norm closure in $L^\infty [0,1]$ of
the linear space of the characteristic functions of open sets $(
\frac{k}{2^n},\frac{k+1}{2^n})$ and of singleton sets $\{
\frac{\ell}{2^n}\}$, $n\geq 0$, $0\leq k<2^n$, $0\leq \ell\leq
2^n$. In particular this shows that $Z(\AA)=C[0,1]$.

The connecting maps $K_0(\AA_n)\hookrightarrow K_0(\AA_{n+1})$
correspond to the polynomial relations
$p_{n+1}(t)=(1+t+t^2)p_n(t^2)$. These polynomials are closely
related to the \emph{Stern--Brocot sequence}. The origins of this
remarkable sequence, which has attracted considerable interest in
time, can be traced back to Eisenstein (see \cite{Ste},
\cite{Bro}, or the contemporary reference \cite{Slo} for a
thorough bibliography on this subject). In our framework the
Stern--Brocot sequence $q(n,k)$, $n\geq 0$, $0\leq k<2^n$, simply
appears as the sizes of the central summands in $\AA_n \cong
\oplus_{k=0}^{2^{n-1}} \M_{q_{(n,k)}}\oplus \C$, where $\M_r$
denotes the $C^*$-algebra of $r\times r$ matrices with complex
entries.

The Bratteli diagram $\GG$ has some apparent symmetries. In the
last section we employ the \emph{AF algebra path model for AF
algebras} to express them, constructing sequences of projections
in $\AA$ that satisfy certain braiding relations reminiscent of
the \emph{Temperley-Lieb-Jones relations}. In particular, for
every $\tau \in ( 0,\frac{1}{4}]$, we construct projections $E_n
\neq 0$ in $\AA$ such that $E_n E_{n\pm 1}E_n \leq \tau E_n$ and
$[E_n,E_m]=0$ if $\vert n-m\vert \geq 2$. This suggests a possible
connection with a class of statistical mechanics models with
partition functions closely related to Riemann's zeta function,
called \emph{Farey spin chains}, that have been studied in recent
years by Knauf, Kleban, and their collaborators (see, e.g.
\cite{Kn,Kn1,Kn2,FK,PFK} and references therein).

\section{The Pascal triangle with memory as a Bratelli diagram}
The {\em Pascal triangle with memory} is a graph $\GG=(\VV,\EE)$
defined as follows:
\begin{itemize}
\item The {\em vertex set} $\VV$ is the disjoint union
$\uplus_{n\geq 0} \VV_n$ of the sets $\VV_n =\{ (n,k) : 0\leq
k\leq 2^n\}$ of {\em vertices at floor $n$}; \item The set of {\em
edges} is defined as $\EE= \uplus_{n\geq 0} \EE_n$, where $\EE_n$
is the set of edges connecting vertices at floor $n$ with those at
floor $n+1$ under the rule that $(n,k)$ is connected with
$(n+1,\ell)$ precisely when $\vert 2k-\ell\vert \leq 1$. There are
no edges connecting vertices from $\VV_i$ and $\VV_j$ when $\vert
i-j\vert \geq 2$.
\end{itemize}

To each vertex $(n,k)$ we attach the {\em label}
$r(n,k)=\frac{p(n,k)}{q(n,k)}$, with non-negative integers
$p(n,k)$, $q(n,k)$ defined recursively for $n\geq 0$ by
\begin{equation*}
\begin{cases}
q(n,0)=q(n,2^n)=1,\quad p(n,0)=0,\quad p(n,2^n)=1; \\
q(n+1,2k)=q(n,k), \quad p(n+1,2k)=p(n,k),\qquad
0\leq k\leq 2^n;\\
q(n+1,2k+1)=q(n,k) +q(n,k+1),\\
p(n+1,2k+1)=p(n,k)+p(n,k+1),\qquad 0\leq k< 2^n.
\end{cases}
\end{equation*}

Note that $r(n,0)=0<r(n,1)=\frac{1}{n+1}<\cdots <r(n,2^n)=1$ gives a
partition of $[0,1]$, and
\begin{equation*}
p(n,k+1)q(n,k)-p(n,k)q(n,k+1)=1, \qquad n\geq 0,\ 0\leq k<2^n,
\end{equation*}
showing in particular that $p(n,k)$ and $q(n,k)$ are relatively
prime.

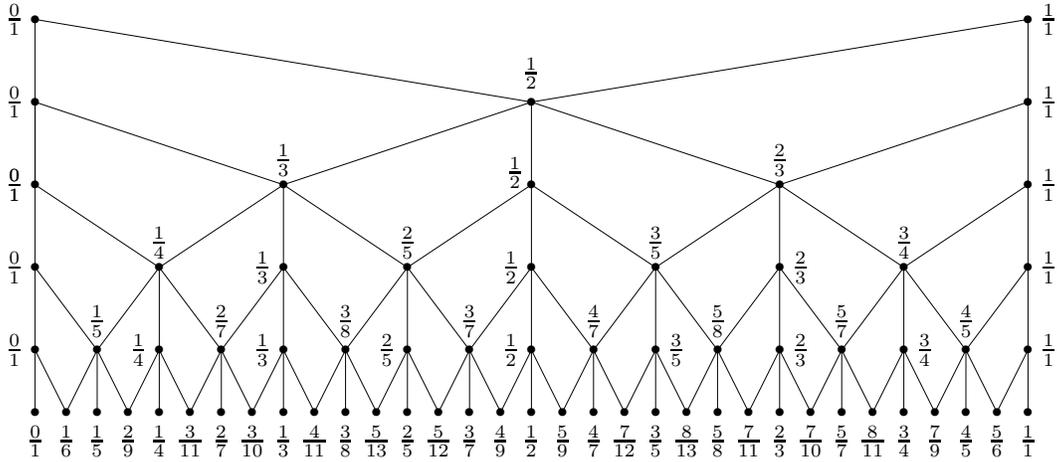
\begin{figure}[htb]
\begin{center}
\unitlength 0.55mm
\begin{picture}(0,110)(0,0)

\path(-120,80)(-120,100)(0,80)(120,100)(120,80)
\path(-120,80)(-120,60)
\path(-120,60)(-120,80)(-60,60)(0,80)(60,60)(120,80)(120,60)
\path(0,80)(0,60)

\path(-120,40)(-120,60)(-90,40)(-60,60)(-30,40)(0,60)(30,40)(60,60)(90,40)(120,60)(120,40)
\path(-60,60)(-60,40) \path(0,60)(0,40) \path(60,60)(60,40)
\path(-120,20)(-120,40)(-105,20)(-90,40)(-75,20)(-60,40)(-45,20)(-30,40)(-15,20)(0,40)(15,20)
(30,40)(45,20)(60,40)(75,20)(90,40)(105,20)(120,40)(120,20)
\path(-90,40)(-90,20) \path(-60,40)(-60,20) \path(-30,40)(-30,20)
\path(0,40)(0,20) \path(30,40)(30,20) \path(60,40)(60,20)
\path(90,40)(90,20)

\path(-120,5)(-120,20)(-112.5,5)(-105,20)(-97.5,5)(-90,20)(-82.5,5)(-75,20)(-67.5,5)(-60,20)
(-52.5,5)(-45,20)(-37.5,5)(-30,20)(-22.5,5)(-15,20)(-7.5,5)(0,20)
\path(120,5)(120,20)(112.5,5)(105,20)(97.5,5)(90,20)(82.5,5)(75,20)(67.5,5)(60,20)
(52.5,5)(45,20)(37.5,5)(30,20)(22.5,5)(15,20)(7.5,5)(0,20)

\path(-105,20)(-105,5) \path(-90,20)(-90,5) \path(-75,20)(-75,5)
\path(-60,20)(-60,5) \path(-45,20)(-45,5) \path(-30,20)(-30,5)
\path(-15,20)(-15,5) \path(0,20)(0,5) \path(105,20)(105,5)
\path(90,20)(90,5) \path(75,20)(75,5) \path(60,20)(60,5)
\path(45,20)(45,5) \path(30,20)(30,5) \path(15,20)(15,5)

\put(-125,100){\makebox(0,0){\small $\frac{0}{1}$}}
\put(0,87){\makebox(0,0){\small $\frac{1}{2}$}}
\put(125,100){\makebox(0,0){\small $\frac{1}{1}$}}
\put(-125,80){\makebox(0,0){\small $\frac{0}{1}$}}
\put(125,80){\makebox(0,0){\small $\frac{1}{1}$}}
\put(-125,60){\makebox(0,0){\small $\frac{0}{1}$}}
\put(-60,66){\makebox(0,0){\small $\frac{1}{3}$}}
\put(-4,63){\makebox(0,0){\small $\frac{1}{2}$}}
\put(60,66){\makebox(0,0){\small $\frac{2}{3}$}}
\put(125,60){\makebox(0,0){\small $\frac{1}{1}$}}
\put(-125,60){\makebox(0,0){\small $\frac{0}{1}$}}

\put(-125,40){\makebox(0,0){\small $\frac{0}{1}$}}
\put(-90,46){\makebox(0,0){\small $\frac{1}{4}$}}
\put(-65,40){\makebox(0,0){\small $\frac{1}{3}$}}
\put(-30,46){\makebox(0,0){\small $\frac{2}{5}$}}
\put(-5,40){\makebox(0,0){\small $\frac{1}{2}$}}
\put(30,46){\makebox(0,0){\small $\frac{3}{5}$}}
\put(65,40){\makebox(0,0){\small $\frac{2}{3}$}}
\put(90,46){\makebox(0,0){\small $\frac{3}{4}$}}
\put(125,40){\makebox(0,0){\small $\frac{1}{1}$}}

\put(-120,100){\makebox(0,0){{\tiny $\bullet$}}}
\put(120,100){\makebox(0,0){{\tiny $\bullet$}}}
\put(-120,80){\makebox(0,0){{\tiny $\bullet$}}}
\put(0,80){\makebox(0,0){{\tiny $\bullet$}}}
\put(120,80){\makebox(0,0){{\tiny $\bullet$}}}
\put(-120,60){\makebox(0,0){{\tiny $\bullet$}}}
\put(-60,60){\makebox(0,0){{\tiny $\bullet$}}}
\put(0,60){\makebox(0,0){{\tiny $\bullet$}}}
\put(60,60){\makebox(0,0){{\tiny $\bullet$}}}
\put(120,60){\makebox(0,0){{\tiny $\bullet$}}}

\put(-120,40){\makebox(0,0){{\tiny $\bullet$}}}
\put(-90,40){\makebox(0,0){{\tiny $\bullet$}}}
\put(-60,40){\makebox(0,0){{\tiny $\bullet$}}}
\put(-30,40){\makebox(0,0){{\tiny $\bullet$}}}
\put(0,40){\makebox(0,0){{\tiny $\bullet$}}}
\put(30,40){\makebox(0,0){{\tiny $\bullet$}}}
\put(60,40){\makebox(0,0){{\tiny $\bullet$}}}
\put(90,40){\makebox(0,0){{\tiny $\bullet$}}}
\put(120,40){\makebox(0,0){{\tiny $\bullet$}}}

\put(-125,20){\makebox(0,0){\small $\frac{0}{1}$}}
\put(-105,27){\makebox(0,0){\small $\frac{1}{5}$}}
\put(-95,20){\makebox(0,0){\small $\frac{1}{4}$}}
\put(-75,27){\makebox(0,0){\small $\frac{2}{7}$}}
\put(-65,20){\makebox(0,0){\small $\frac{1}{3}$}}
\put(-45,27){\makebox(0,0){\small $\frac{3}{8}$}}
\put(-35,20){\makebox(0,0){\small $\frac{2}{5}$}}
\put(-15,27){\makebox(0,0){\small $\frac{3}{7}$}}
\put(-5,20){\makebox(0,0){\small $\frac{1}{2}$}}
\put(125,20){\makebox(0,0){\small $\frac{1}{1}$}}
\put(105,27){\makebox(0,0){\small $\frac{4}{5}$}}
\put(95,20){\makebox(0,0){\small $\frac{3}{4}$}}
\put(75,27){\makebox(0,0){\small $\frac{5}{7}$}}
\put(65,20){\makebox(0,0){\small $\frac{2}{3}$}}
\put(45,27){\makebox(0,0){\small $\frac{5}{8}$}}
\put(35,20){\makebox(0,0){\small $\frac{3}{5}$}}
\put(15,27){\makebox(0,0){\small $\frac{4}{7}$}}

\put(-120,20){\makebox(0,0){{\tiny $\bullet$}}}
\put(-105,20){\makebox(0,0){{\tiny $\bullet$}}}
\put(-90,20){\makebox(0,0){{\tiny $\bullet$}}}
\put(-75,20){\makebox(0,0){{\tiny $\bullet$}}}
\put(-60,20){\makebox(0,0){{\tiny $\bullet$}}}
\put(-45,20){\makebox(0,0){{\tiny $\bullet$}}}
\put(-30,20){\makebox(0,0){{\tiny $\bullet$}}}
\put(-15,20){\makebox(0,0){{\tiny $\bullet$}}}
\put(0,20){\makebox(0,0){{\tiny $\bullet$}}}
\put(120,20){\makebox(0,0){{\tiny $\bullet$}}}
\put(105,20){\makebox(0,0){{\tiny $\bullet$}}}
\put(90,20){\makebox(0,0){{\tiny $\bullet$}}}
\put(75,20){\makebox(0,0){{\tiny $\bullet$}}}
\put(60,20){\makebox(0,0){{\tiny $\bullet$}}}
\put(45,20){\makebox(0,0){{\tiny $\bullet$}}}
\put(30,20){\makebox(0,0){{\tiny $\bullet$}}}
\put(15,20){\makebox(0,0){{\tiny $\bullet$}}}

\put(-120,-2){\makebox(0,0){\small $\frac{0}{1}$}}
\put(-112.5,-2){\makebox(0,0){\small $\frac{1}{6}$}}
\put(-105,-2){\makebox(0,0){\small $\frac{1}{5}$}}
\put(-97.5,-2){\makebox(0,0){\small $\frac{2}{9}$}}
\put(-90,-2){\makebox(0,0){\small $\frac{1}{4}$}}
\put(-82.5,-2){\makebox(0,0){\small $\frac{3}{11}$}}
\put(-75,-2){\makebox(0,0){\small $\frac{2}{7}$}}
\put(-67.5,-2){\makebox(0,0){\small $\frac{3}{10}$}}
\put(-60,-2){\makebox(0,0){\small $\frac{1}{3}$}}
\put(-52.5,-2){\makebox(0,0){\small $\frac{4}{11}$}}
\put(-45,-2){\makebox(0,0){\small $\frac{3}{8}$}}
\put(-37.5,-2){\makebox(0,0){\small $\frac{5}{13}$}}
\put(-30,-2){\makebox(0,0){\small $\frac{2}{5}$}}
\put(-22.5,-2){\makebox(0,0){\small $\frac{5}{12}$}}
\put(-15,-2){\makebox(0,0){\small $\frac{3}{7}$}}
\put(-7.5,-2){\makebox(0,0){\small $\frac{4}{9}$}}
\put(0,-2){\makebox(0,0){\small $\frac{1}{2}$}}
\put(120,-2){\makebox(0,0){\small $\frac{1}{1}$}}
\put(112.5,-2){\makebox(0,0){\small $\frac{5}{6}$}}
\put(105,-2){\makebox(0,0){\small $\frac{4}{5}$}}
\put(97.5,-2){\makebox(0,0){\small $\frac{7}{9}$}}
\put(90,-2){\makebox(0,0){\small $\frac{3}{4}$}}
\put(82.5,-2){\makebox(0,0){\small $\frac{8}{11}$}}
\put(75,-2){\makebox(0,0){\small $\frac{5}{7}$}}
\put(67.5,-2){\makebox(0,0){\small $\frac{7}{10}$}}
\put(60,-2){\makebox(0,0){\small $\frac{2}{3}$}}
\put(52.5,-2){\makebox(0,0){\small $\frac{7}{11}$}}
\put(45,-2){\makebox(0,0){\small $\frac{5}{8}$}}
\put(37.5,-2){\makebox(0,0){\small $\frac{8}{13}$}}
\put(30,-2){\makebox(0,0){\small $\frac{3}{5}$}}
\put(22.5,-2){\makebox(0,0){\small $\frac{7}{12}$}}
\put(15,-2){\makebox(0,0){\small $\frac{4}{7}$}}
\put(7.5,-2){\makebox(0,0){\small $\frac{5}{9}$}}

\put(-120,5){\makebox(0,0){{\tiny $\bullet$}}}
\put(-112.5,5){\makebox(0,0){{\tiny $\bullet$}}}
\put(-105,5){\makebox(0,0){{\tiny $\bullet$}}}
\put(-97.5,5){\makebox(0,0){{\tiny $\bullet$}}}
\put(-90,5){\makebox(0,0){{\tiny $\bullet$}}}
\put(-82.5,5){\makebox(0,0){{\tiny $\bullet$}}}
\put(-75,5){\makebox(0,0){{\tiny $\bullet$}}}
\put(-67.5,5){\makebox(0,0){{\tiny $\bullet$}}}
\put(-60,5){\makebox(0,0){{\tiny $\bullet$}}}
\put(-52.5,5){\makebox(0,0){{\tiny $\bullet$}}}
\put(-45,5){\makebox(0,0){{\tiny $\bullet$}}}
\put(-37.5,5){\makebox(0,0){{\tiny $\bullet$}}}
\put(-30,5){\makebox(0,0){{\tiny $\bullet$}}}
\put(-22.5,5){\makebox(0,0){{\tiny $\bullet$}}}
\put(-15,5){\makebox(0,0){{\tiny $\bullet$}}}
\put(-7.5,5){\makebox(0,0){{\tiny $\bullet$}}}
\put(0,5){\makebox(0,0){{\tiny $\bullet$}}}
\put(120,5){\makebox(0,0){{\tiny $\bullet$}}}
\put(112.5,5){\makebox(0,0){{\tiny $\bullet$}}}
\put(105,5){\makebox(0,0){{\tiny $\bullet$}}}
\put(97.5,5){\makebox(0,0){{\tiny $\bullet$}}}
\put(90,5){\makebox(0,0){{\tiny $\bullet$}}}
\put(82.5,5){\makebox(0,0){{\tiny $\bullet$}}}
\put(75,5){\makebox(0,0){{\tiny $\bullet$}}}
\put(67.5,5){\makebox(0,0){{\tiny $\bullet$}}}
\put(60,5){\makebox(0,0){{\tiny $\bullet$}}}
\put(52.5,5){\makebox(0,0){{\tiny $\bullet$}}}
\put(45,5){\makebox(0,0){{\tiny $\bullet$}}}
\put(37.5,5){\makebox(0,0){{\tiny $\bullet$}}}
\put(30,5){\makebox(0,0){{\tiny $\bullet$}}}
\put(22.5,5){\makebox(0,0){{\tiny $\bullet$}}}
\put(15,5){\makebox(0,0){{\tiny $\bullet$}}}
\put(7.5,5){\makebox(0,0){{\tiny $\bullet$}}}

\end{picture}
\end{center}
\caption{The Pascal triangle with memory $\GG$} \label{Figure2}
\end{figure}

Conversely, for every pair $\frac{p}{q}<\frac{p'}{q'}$ of rational
numbers with $p'q-pq'=1$, $0\leq p\leq q$ and $0\leq p'\leq q'$,
there exists a unique pair of integers $(n,k)$ with $n\geq 0$,
$0\leq k<2^n$, such that $r(n,k)=\frac{p}{q}$ and
$r(n,k+1)=\frac{p'}{q'}$. This correspondence establishes a
bijection between the vertices from $\VV\setminus \{ (n,2^n):n\geq
0\}$ and the set
\begin{equation*}
\Gamma^+=\left\{ \left[ \begin{smallmatrix} p' & p \\ q' & q
\end{smallmatrix}\right] \in SL_2 (\Z):0\leq p\leq q,\ 0\leq
p'\leq q'\right\} \subset SL_2(\Z).
\end{equation*}

\begin{remark}\label{R?}
{\em The mapping $r(n,k)\mapsto \frac{k}{2^n}$, $0\leq k\leq 2^n$,
$n\geq 0$, extends by continuity to \emph{Minkowski's question
mark function} $\mbox{\large ?}:[0,1]\rightarrow [0,1]$ defined on
(reduced) continued fractions as
\begin{equation*}
\mbox{\large ?} ([a_1,a_2,\ldots]) \ =\ \sum\limits_{k\geq 1}
\frac{(-1)^{k-1}}{2^{(a_1+\cdots +a_k)-1}}.
\end{equation*}
The map $\mbox{\large ?}$ is strictly increasing and singular, and
establishes remarkable one-to-one correspondences between rational
and dyadic numbers, and respectively between quadratic irrationals
and rational numbers in $[0,1]$ (see \cite{Min,Den,Sal}). }
\end{remark}

In this paper we shall consider the AF algebra $\AA$ associated
with the Bratteli diagram $D (\AA)=\GG$ from Figure \ref{Figure2}.
For the connection between Bratteli diagrams, AF algebras, and
their ideals, we refer to the classical reference \cite{Bra}. We
write $(n,k)\downarrow (n^\prime,k^\prime)$ when $n^\prime =n+1$
and there is at least one edge between the vertices $(n,k)$ and
$(n^\prime,k^\prime)$ in the Bratteli diagram, and
$(n,k)\Downarrow (n^\prime,k^\prime)$ when $n<n^\prime$ and there
are vertices
$(n,k_0=k),(n+1,k_1),\ldots,(n^\prime,k_{n^\prime-n}=k^\prime)$
such that $(n+r,k_r)\downarrow (n+r+1,k_{r+1})$,
$r=0,\ldots,n^\prime-n-1$. In algebraic terms this is equivalent
to $e_{(n,k)} e_{(n^\prime,k^\prime)} \neq 0$, where $e_{(n,k)}$
denotes the central projection in $\AA_n$ that corresponds to the
vertex $(n,k)$ of the diagram. The AF algebra $\AA$ is the
inductive limit $\varinjlim \AA_n$, where
\begin{equation*}
\AA_n =\mbox{\small $\displaystyle \bigoplus\limits_{0\leq k\leq
2^n}$} \M_{q(n,k)}
\end{equation*}
and each embedding  $\AA_n \hookrightarrow \AA_{n+1}$ is given by
the Bratteli diagram from Figure \ref{Figure2}.

\begin{remark}\label{R22}
{\em Consider the set $\VV_\ast$ of vertices of $\GG$ of form
$(n,k)$ with $0\leq k\leq 2^n$ and $k$ odd, and the map
$\Phi:\VV_\ast\rightarrow \N$, $\Phi(n,k)=q(n,k)$. The inverse
image $\Phi^{-1}(q)$ of $q$ contains exactly $\varphi (q)$
elements, where $\varphi$ denotes Euler's totient function; in
particular $q$ is prime if and only if $\# \Phi^{-1}(q)=q-1$. This
remark shows, cf. \cite{Kn}, that the partition function
associated with the corresponding Farey spin chain is
$\sum_{n=1}^\infty \varphi(n)n^{-s}$, which is equal to $\zeta
(s-1)/\zeta(s)$ when $\Re s>2$.}
\end{remark}

\begin{remark}\label{R1}
{\em (i) The integers $q(n,k)$ satisfy the equality}
\begin{equation*}
\sum\limits_{0\leq k\leq 2^n} q(n,k)=3^n +1.
\end{equation*}

{\em (ii) Consider the Bratteli diagram obtained by deleting in
$\GG$ all vertices $(n,0)$ and denote the corresponding AF algebra
by $\BB=\varinjlim \BB_n$. It is clear that $\BB$ is an ideal in
$\AA$ and $\AA /\BB \cong \C$. Moreover, }
\begin{equation*}
\BB_n=\mbox{\small $\displaystyle \bigoplus\limits_{1\leq k\leq
2^n}$} \M_{p(n,k)} ,
\end{equation*}
{\em thus the ranks of the central summands of the building blocks
of $\BB$ give the complete list of numerators $p(n,k)$. We also
have}
\begin{equation*}
\sum\limits_{0\leq k\leq 2^n} p(n,k)=\frac{3^n+1}{2}  .
\end{equation*}
\end{remark}

\section{The primitive ideal space of the AF algebra $\AA$}

We denote $ \I=\{ \theta \in (0,1):\theta \notin \Q\}$ and
$\Q_{(0,1)}= \Q \cap (0,1)$.

\begin{figure}[htb]
\begin{center}
\unitlength 0.6mm
\begin{picture}(20,40)(0,7)

\path(-43,10)(-17,10) \path(-43,9)(-17,9) \path(-43,11)(-17,11)
\path(-44,11)(-16,39) \path(-13,10.5)(13,10.5)
\path(-13,9.5)(13,9.5) \path(17,10)(43,10) \path(17,9)(43,9)
\path(17,11)(43,11) \path(46,10)(74,10) \path(-14,39)(14,11)
\path(16,39)(44,11) \path(46,39)(74,11) \path(-14,11)(14,39)
\path(16,11)(44,39) \path(46,11)(74,39)

\put(-45,10){\makebox(0,0){{\footnotesize $\bullet$}}}
\put(-15,10){\makebox(0,0){{\footnotesize $\bullet$}}}
\put(-15,40){\makebox(0,0){{\footnotesize $\bullet$}}}
\put(15,10){\makebox(0,0){{\footnotesize $\bullet$}}}
\put(45,10){\makebox(0,0){{\footnotesize $\bullet$}}}
\put(75,10){\makebox(0,0){{\footnotesize $\bullet$}}}
\put(15,40){\makebox(0,0){{\footnotesize $\bullet$}}}
\put(45,40){\makebox(0,0){{\footnotesize $\bullet$}}}
\put(75,40){\makebox(0,0){{\footnotesize $\bullet$}}}
\put(-30,5){\makebox(0,0){{$a_1$}}}
\put(0,5){\makebox(0,0){{$a_2$}}} \put(30,5){\makebox(0,0){{
$a_3$}}} \put(60,5){\makebox(0,0){{$a_4$}}}
\put(95,10){\makebox(0,0){{$\ldots$}}}
\put(95,40){\makebox(0,0){{$\ldots$}}}
\end{picture}
\end{center}
\caption{The Bratteli diagram $D(\CC_\theta)$} \label{Figure3}
\end{figure}
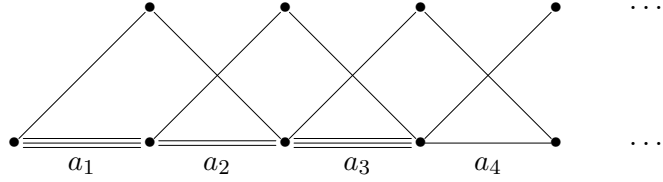

The $C^*$-algebra $\AA$ is not simple and has a rich (and
potentially interesting) structure of ideals. We first relate
$\AA$ with the AF algebra $\CC_\theta$ associated by Effros and
Shen \cite{ES} to the continued fraction decomposition $\theta
=[a_1,a_2,\ldots]$ of $\theta \in \I$. The Bratteli diagram
$D(\CC_\theta)$ of the simple $C^*$-algebra $\CC_\theta$ is given
in Figure \ref{Figure3}.

\begin{figure}[htb]
\begin{center}
\unitlength 0.6mm
\begin{picture}(20,40)(0,7)

\texture{c 00000} \shade\path(-50,5)(-50,15)(90,15)(90,5)(-50,5)
\path(-45,10)(85,10) \path(-45,40)(85,40) \path(-45,40)(-15,10)
\path(-15,40)(15,10) \path(15,40)(45,10) \path(45,40)(75,10)
\path(75,40)(85,30)

\put(-45,40){\makebox(0,0){{\footnotesize
$\bullet$}}}\put(-45,10){\makebox(0,0){{\footnotesize $\bullet$}}}
\put(-15,10){\makebox(0,0){{\footnotesize $\bullet$}}}
\put(-15,40){\makebox(0,0){{\footnotesize $\bullet$}}}
\put(15,10){\makebox(0,0){{\footnotesize $\bullet$}}}
\put(45,10){\makebox(0,0){{\footnotesize $\bullet$}}}
\put(75,10){\makebox(0,0){{\footnotesize $\bullet$}}}
\put(15,40){\makebox(0,0){{\footnotesize $\bullet$}}}
\put(45,40){\makebox(0,0){{\footnotesize $\bullet$}}}
\put(75,40){\makebox(0,0){{\footnotesize $\bullet$}}}

\put(-55,40){\makebox(0,0){{$\C$}}}
\put(-55,10){\makebox(0,0){{$\C$}}}

\end{picture}
\end{center}

\caption{The Bratteli diagram of the $C^*$-algebra of unitized
compact operators} \label{Figure31}
\end{figure}
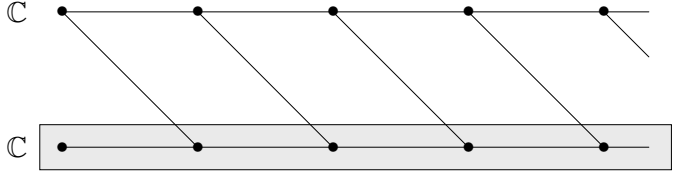

The $C^*$-algebra of unitized compact operators $\widetilde{\K}=\C
I+\K$ is an AF algebra and we have a short exact sequence
$0\rightarrow \K\rightarrow \widetilde{\K}\rightarrow \C
\rightarrow 0$, made explicit by the Bratteli diagram in Figure
\ref{Figure31}, where the shaded subdiagram corresponds to the
ideal $\K$. Replacing $\C\oplus \C$ by $\M_q\oplus \M_{q^\prime}$
one gets an AF algebra $\AA_{(q,q^\prime)}$ which is an extension
of $\K$ by $\M_q$.

We first show that Effros-Shen algebras arise naturally as
quotients of our AF algebra $\AA$ and that the corresponding
ideals belong to the primitive ideal space $\Prim \AA$. The
\emph{Farey map} $F:[0,1]\rightarrow [0,1]$ defined \cite{Fe} by
\begin{equation}\label{2.1}
F(x)=\begin{cases} \vspace{.2cm} \displaystyle \frac{x}{1-x} &
\mbox{\rm if $x\in \big[ 0,\frac{1}{2}\big],$} \\ \displaystyle
\frac{1-x}{x} & \mbox{\rm if $x\in \big( \frac{1}{2},1\big]$,}
\end{cases}
\end{equation}
acts on infinite (reduced) continued fractions as
\begin{equation*}
F([a_1,a_2,a_3,\ldots])=[a_1-1,a_2,a_3,\ldots ].
\end{equation*}
For each $y\in [0,1]$ the equation $F(x)=y$ has exactly two
solutions $x\in [0,1]$ given by
\begin{equation}\label{2.2}
x=F_1(y)=\frac{y}{1+y}\quad \mbox{\rm and}\quad
x=F_2(y)=\frac{1}{1+y}=1-F_1(y).
\end{equation}
One has $F_1([a_1,a_2,\ldots])=[a_1+1,a_2,\ldots]$ and
$F_2([a_1,a_2,\ldots])=[1,a_1,a_2,\ldots]$. Rational numbers are
generated by the backwards orbit of $F$ as follows:
\begin{equation*}
\{ F^{-n}(\{0\}):n=0,1,2,\ldots\}=\Q \cap [0,1].
\end{equation*}
More precisely, for each $n\in\N$ one has
\begin{equation*}
\begin{split}
F^{-n}(\{ 0\}) & = \{ r(n-1,k):0\leq k\leq 2^{n-1}\} \\
& =\big\{F_{i_1}^{\alpha_1} \ldots F_{i_k}^{\alpha_k} (0):i_j\in\{
1,2\},\ i_1\neq \cdots\neq i_k,\ \alpha_1+\cdots+\alpha_k=n\big\}
\\ & =\big\{ [a_1,\ldots,a_r]: a_1+\cdots+a_r \leq n\big\}.
\end{split}
\end{equation*}

In the next statement, given relatively prime integers $0<p<q$,
$\overline{p}$ will denote the multiplicative inverse of $p$
modulo $q$ in $\{ 1,\ldots,q-1\}$.

\begin{proposition}\label{P2}
{\em (i)} For each $\theta \in \I$, there is $I_\theta \in \Prim
\AA$ such that $\ \AA/I_\theta \cong \CC_{\theta}$.

{\em (ii)} Given $\theta=\frac{p}{q}\in \Q_{(0,1)}$ in lowest
terms, there are $I_\theta, I_\theta^+,I_\theta^-\in\Prim\AA$ such
that $\AA/I_\theta \cong \M_q$, $\AA/I^-_\theta\cong
\AA_{(q,\overline{p})}$, and $\AA/I^+_\theta \cong
\AA_{(q,q-\overline{p})}$.

{\em (iii)} There are $I_0, I_0^+, I_1, I_1^-\in \Prim\AA$ such
that $\AA/I_0 \cong \AA/I_1 \cong \C$ and $\AA/I_0^+ \cong
\AA/I_1^- \cong \widetilde{\K}$.
\end{proposition}

\begin{proof} (i) Let $\theta \in \I$ with continued fraction
$[a_1,a_2,\ldots]$ and $r_\ell=r_\ell (\theta)=p_\ell /q_\ell
=[a_1,\ldots,a_\ell ]$ be its $\ell^{\mathrm{th}}$ convergent,
where $p_\ell =p_\ell (\theta)$ and $q_\ell =q_\ell (\theta)$ can
be recursively defined by
\begin{equation*}
\begin{cases} \vspace{0.2cm}
p_{-1}=1,\ q_{-1}=0,\quad p_0=0,\ q_0=1; \\
\left[ \begin{matrix} p_\ell & q_\ell \\ p_{\ell -1} & q_{\ell -1}
\end{matrix} \right] =\left[ \begin{matrix} a_\ell & 1 \\ 1 & 0
\end{matrix} \right] \left[ \begin{matrix} p_{\ell -1} & q_{\ell -1} \\
p_{\ell -2} & q_{\ell -2} \end{matrix} \right] ,\qquad \ell\geq 1.
\end{cases}
\end{equation*}
The relation $p_\ell q_{\ell -1}-p_{\ell -1} q_\ell =(-1)^{\ell
-1}$ shows in particular that $\gcd (p_\ell ,q_\ell )=1$.

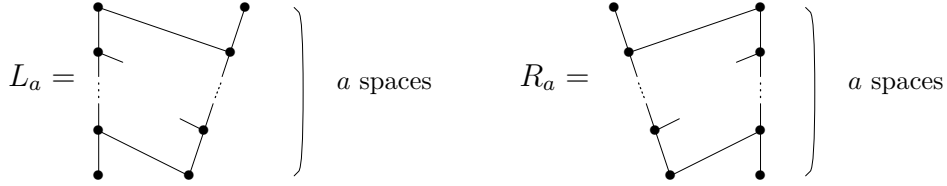
\begin{figure}[htb]
\begin{center}
\unitlength 0.4mm
\begin{picture}(120,50)(10,20)

\path(-45,30)(-15,15) \path(-10,30)(-18,34) \path(-45,71)(-1.3,56)
\path(-45,56)(-37,52.75) \path(-45,15)(-45,39)
\path(-45,48)(-45,71) \path(-15,15)(-7,39)

\dottedline{2}(-45,41)(-45,46) \path(-4,48)(3.66,71)
\dottedline{2}(-6.3,41)(-4.4,46)

\put(-75,45){\large $L_a=$} \put(50,45){\makebox(0,0){{\small $a$
spaces}}} \put(-45,15){\makebox(0,0){{\footnotesize $\bullet$}}}
\put(-45,30){\makebox(0,0){{\footnotesize $\bullet$}}}
\put(-45,56){\makebox(0,0){{\footnotesize $\bullet$}}}
\put(-45,71){\makebox(0,0){{\footnotesize $\bullet$}}}
\put(-15,15){\makebox(0,0){{\footnotesize $\bullet$}}}
\put(-10,30){\makebox(0,0){{\footnotesize $\bullet$}}}
\put(-1.3,56){\makebox(0,0){{\footnotesize $\bullet$}}}
\put(3.66,71){\makebox(0,0){{\footnotesize $\bullet$}}}

\spline(20,15)(23,18)(23,68)(20,71)

\path(175,30)(145,15) \path(140,30)(148,34)
\path(175,71)(131.3,56) \path(175,56)(167,52.75)
\path(175,15)(175,39) \path(175,48)(175,71) \path(145,15)(137,39)

\dottedline{2}(175,41)(175,46) \path(134,48)(126.33,71)
\dottedline{2}(136.3,41)(134.4,46)

\put(95,45){\large $R_a=$}

\put(175,15){\makebox(0,0){{\footnotesize $\bullet$}}}
\put(175,30){\makebox(0,0){{\footnotesize $\bullet$}}}
\put(175,56){\makebox(0,0){{\footnotesize $\bullet$}}}
\put(175,71){\makebox(0,0){{\footnotesize $\bullet$}}}
\put(145,15){\makebox(0,0){{\footnotesize $\bullet$}}}
\put(140,30){\makebox(0,0){{\footnotesize $\bullet$}}}
\put(131.3,56){\makebox(0,0){{\footnotesize $\bullet$}}}
\put(126.33,71){\makebox(0,0){{\footnotesize $\bullet$}}}

\spline(190,15)(193,18)(193,68)(190,71)
\put(220,45){\makebox(0,0){{\small $a$ spaces}}}

\end{picture}
\end{center}
\caption{The diagrams $L_a$ and $R_a$} \label{Figure4}
\end{figure}

For each $a\in \N=\{ 1,2,\ldots \}$ consider the diagrams $L_a$
and $R_a$ from Figure \ref{Figure4}. Also set $L_0=R_0=\emptyset$.
Clearly $L_{a+b}$ coincides with the concatenation $L_a\circ L_b$
of $L_a$ followed by $L_b$, and we also have $R_{a+b}=R_a\circ
R_b$. Using the obvious identifications between $L_a \circ R_b$,
$R_a\circ L_b$ and $C_a \circ C_b$ (see Figure \ref{Figure5}), and
\eqref{2.2}, we see that the AF algebras generated by
$L_{a_1}\circ R_{a_2}\circ L_{a_3}\circ R_{a_4}\circ\dots$ and
$R_{a_1}\circ L_{a_2}\circ R_{a_3}\circ L_{a_4}\circ\dots$ are
isomorphic to $\CC_{[a_1+1,a_2,a_3,\ldots]}\simeq
\CC_{F_1(\theta)}\simeq \CC_{F_2(\theta)}\simeq
\CC_{[1,a_1,a_2,\ldots]}$ (note that the AF algebra defined by
$C_{a_1}\circ C_{a_2}\circ C_{a_3} \circ \cdots$ is isomorphic to
$\CC_{[a_1+1,a_2,a_3,\ldots]}$).

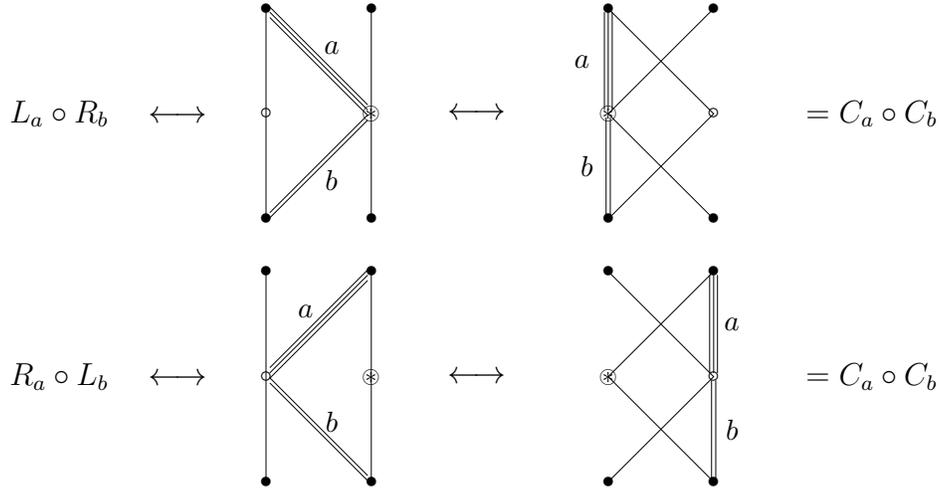
\begin{figure}[htb]
\begin{center}
\unitlength 0.35mm
\begin{picture}(150,170)(0,-80)
%\thinlines

\put(-60,50){\makebox(0,0){\large $L_a\circ R_b\quad
\longleftrightarrow $}}

\path(0,90)(0,10) \path(40,90)(40,10) \path(1.5,88.5)(38.5,51.5)
\path(1.5,90.5)(38.5,53.5) \path(1.5,86.5)(38.5,49.5)
\path(1.5,10.5)(38.5,47.5) \path(1.5,12.5)(38.5,49.5)

\put(40,90){\makebox(0,0){{\footnotesize $\bullet$}}}
\put(40,50){\makebox(0,0){{\footnotesize $\circledast$}}}
\put(40,10){\makebox(0,0){{\footnotesize $\bullet$}}}
\put(0,10){\makebox(0,0){{\footnotesize $\bullet$}}}
\put(0,90){\makebox(0,0){{\footnotesize $\bullet$}}}
\put(0,50){\makebox(0,0){{\footnotesize $\circ$}}}

\put(25,75){\makebox(0,0){$a$}} \put(23,25){\makebox(0,0){ $b$}}

\put(80,50){\makebox(0,0){\large $\longleftrightarrow$}}

\path(170,90)(130,50)(170,10) \path(130,90)(170,50)(130,10)
\path(130,88.5)(130,51.5) \path(128.5,88.5)(128.5,51.5)
\path(131.5,88.5)(131.5,51.5)

\path(129.2,48.5)(129.2,11.5) \path(130.8,48.5)(130.8,11.5)

\put(130,90){\makebox(0,0){{\footnotesize $\bullet$}}}
\put(170,50){\makebox(0,0){{\footnotesize $\circ$}}}
\put(130,50){\makebox(0,0){{\footnotesize $\circledast$}}}
\put(170,90){\makebox(0,0){{\footnotesize $\bullet$}}}
\put(130,10){\makebox(0,0){{\footnotesize $\bullet$}}}
\put(170,10){\makebox(0,0){{\footnotesize $\bullet$}}}
\put(120,70){\makebox(0,0){$a$}} \put(120,30){\makebox(0,0){ $b$}}
\put(230,50){\makebox(0,0){\large $=C_a\circ C_b$}}

\put(-60,-50){\makebox(0,0){\large $R_a\circ L_b\quad
\longleftrightarrow $}}

\path(0,-10)(0,-90) \path(40,-10)(40,-90)

\path(1.5,-48.5)(38.5,-11.5) \path(1.5,-46.5)(38.5,-9.5)
\path(1.5,-50.5)(38.5,-13.5) \path(1.5,-50.5)(38.5,-87.5)
\path(1.5,-52.5)(38.5,-89.5)

\put(40,-10){\makebox(0,0){{\footnotesize $\bullet$}}}
\put(40,-50){\makebox(0,0){{\footnotesize $\circledast$}}}
\put(40,-90){\makebox(0,0){{\footnotesize $\bullet$}}}
\put(0,-90){\makebox(0,0){{\footnotesize $\bullet$}}}
\put(0,-10){\makebox(0,0){{\footnotesize $\bullet$}}}
\put(0,-50){\makebox(0,0){{\footnotesize $\circ$}}}

\put(15,-25){\makebox(0,0){$a$}} \put(25,-67){\makebox(0,0){$b$}}

\put(80,-50){\makebox(0,0){\large $\longleftrightarrow$}}

\path(170,-10)(130,-50)(170,-90) \path(130,-10)(170,-50)(130,-90)
\path(170,-11.5)(170,-48.5) \path(168.5,-11.5)(168.5,-48.5)
\path(171.5,-11.5)(171.5,-48.5)

\path(169.2,-51.5)(169.2,-88.5) \path(170.8,-51.5)(170.8,-88.5)

\put(130,-10){\makebox(0,0){{\footnotesize $\bullet$}}}
\put(170,-50){\makebox(0,0){{\footnotesize $\circ$}}}
\put(130,-50){\makebox(0,0){{\footnotesize $\circledast$}}}
\put(170,-10){\makebox(0,0){{\footnotesize $\bullet$}}}
\put(130,-90){\makebox(0,0){{\footnotesize $\bullet$}}}
\put(170,-90){\makebox(0,0){{\footnotesize $\bullet$}}}
\put(177,-30){\makebox(0,0){{$a$}}}
\put(177,-70){\makebox(0,0){{$b$}}}
\put(230,-50){\makebox(0,0){\large $=C_a\circ C_b$}}
\end{picture}
\end{center}
\caption{The identification between $L_a\circ R_b$, $R_a\circ
L_b$, and $C_a\circ C_b$} \label{Figure5}
\end{figure}

The Bratteli subdiagram $\GG_\theta$ of $\GG$ containing the
vertices $(0,0)$ and $(0,1)$ and defined by $L_{a_1-1}\circ
R_{a_2}\circ L_{a_3}\circ R_{a_4}\circ \cdots$ generates a copy of
$\CC_\theta$. The complement $\GG\setminus \GG_\theta$ is a
directed and hereditary Bratteli diagram as in \cite[Lemma
3.2]{Bra} (see also Figure \ref{Figure7}). Thus there is an ideal
$I_\theta$ in $\AA$ such that $D(I_\theta)=\GG\setminus
\GG_\theta$, $D(\AA/I_\theta)=\GG_\theta$, and $\AA /I_\theta
\cong \CC_{\theta}$. Moreover $I_\theta$ is a primitive ideal cf.
\cite[Theorem 3.8]{Bra}.

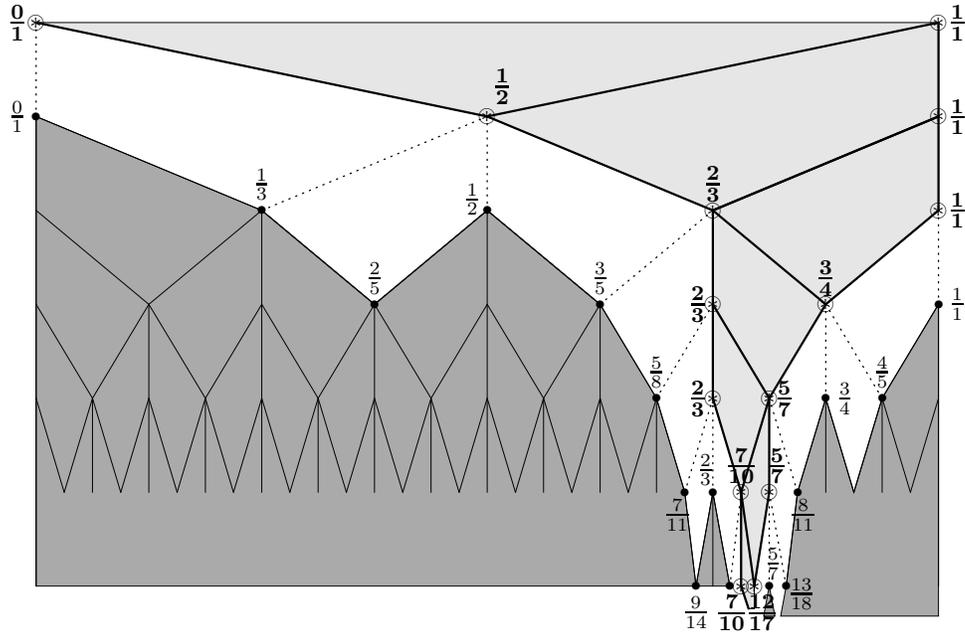
\begin{figure}[htb]
\begin{center}
\unitlength 0.5mm
\begin{picture}(0,160)(5,-25)
%\thinlines
\linethickness{10mm} \texture{cc 0}
\shade\path(-120,-20)(-120,105)(-60,80)(-30,55)(0,80)(30,55)(45,30)(52.5,5)(55.5,-20)(-120,-20)
\shade\path(55.5,-20)(64.5,-20)(60,5)(55.5,-20)
\shade\path(120,-28)(120,55)(105,30)(97.5,5)(90,30)(82.5,5)(79.5,-20)(78,-28)(120,-28)
\shade\path(75,-20)(73.5,-28)(76.5,-28)(75,-20)

\texture{c 0000}
\shade\path(-120,130)(120,130)(120,80)(90,55)(75,30)(75,5)(71,-20)(71,-26)(69.5,-26)(67.5,-20)(67.5,5)(60,30)(60,80)(0,105)(-120,130)

\drawline(60,5)(60,-20)
\path(-120,80)(-90,55)(-60,80)(-30,55)(0,80)(30,55)

\thicklines
\drawline(-120,130)(0,105)(120,130)(120,105)(60,80)(0,105)
\drawline(120,105)(60,80)(60,30)
\drawline(120,105)(120,80)(90,55)(60,80)
\drawline(67.5,-20)(67.5,5)
\drawline(60,30)(67.5,5)(71,-20)(75,5)(75,30)(60,55)
\drawline(67.5,5)(75,30)(90,55) \drawline(67.5,-20)(69.5,-26)

\thinlines \dottedline{2}(-120,105)(-120,130)
\dottedline{2}(-60,80)(0,105)(0,80) \dottedline{2}(30,55)(60,80)
\dottedline{2}(120,55)(120,80) \dottedline{2}(45,30)(60,55)
\dottedline{2}(90,30)(90,55)(105,30)
\dottedline{2}(67.5,5)(64.5,-20)
\dottedline{2}(75,-20)(75,5)(79.5,-20)
\dottedline{2}(60,5)(60,30)(52.5,5) \dottedline{2}(75,30)(82.5,5)
\drawline(90,30)(90,5) \drawline(105,5)(105,30)(112.5,5)(120,30)
\drawline(45,30)(45,5) \drawline(30,55)(30,5)
\drawline(15,30)(15,5) \drawline(0,80)(0,5)
\drawline(-15,30)(-15,5) \drawline(-30,55)(-30,5)
\drawline(-45,30)(-45,5) \drawline(-60,80)(-60,5)
\drawline(-75,30)(-75,5) \drawline(-90,55)(-90,5)
\drawline(-105,30)(-105,5)
\drawline(-120,30)(-112.5,5)(-105,30)(-97.5,5)(-90,30)(-82.5,5)(-75,30)(-67.5,5)(-60,30)
\drawline(-60,30)(-52.5,5)(-45,30)(-37.5,5)(-30,30)(-22.5,5)(-15,30)(-7.5,5)(0,30)(7.5,5)(15,30)
\drawline(15,30)(22.5,5)(30,30)(37.5,5)(45,30)
\drawline(-120,55)(-105,30)(-90,55)(-75,30)(-60,55)(-45,30)(-30,55)(-15,30)(0,55)(15,30)(30,55)

\drawline(-120,-20)(-120,105)(-60,80)(-30,55)(0,80)(30,55)(45,30)(52.5,5)(55.5,-20)(60,5)(64.5,-20)
\drawline(120,-20)(120,55)(105,30)(97.5,5)(90,30)(82.5,5)(79.5,-20)

\put(-125,130){\makebox(0,0){{$\mathbf{\frac{0}{1}}$}}}
\put(4,113){\makebox(0,0){{$\mathbf{\frac{1}{2}}$}}}
\put(125,130){\makebox(0,0){{$\mathbf{\frac{1}{1}}$}}}
\put(-125,105){\makebox(0,0){{\small $\frac{0}{1}$}}}
\put(125,105){\makebox(0,0){{$\mathbf{\frac{1}{1}}$}}}
\put(-60,87){\makebox(0,0){{\small $\frac{1}{3}$}}}
\put(-4,83){\makebox(0,0){{\small $\frac{1}{2}$}}}
\put(60,87){\makebox(0,0){{$\mathbf{\frac{2}{3}}$}}}
\put(125,80){\makebox(0,0){{$\mathbf{\frac{1}{1}}$}}}

\put(-30,62){\makebox(0,0){{\small $\frac{2}{5}$}}}
\put(30,62){\makebox(0,0){{\small $\frac{3}{5}$}}}
\put(56,55){\makebox(0,0){{$\mathbf{\frac{2}{3}}$}}}
\put(90,62){\makebox(0,0){{$\mathbf{\frac{3}{4}}$}}}
\put(125,55){\makebox(0,0){{\small $\frac{1}{1}$}}}

\put(-120,130){\makebox(0,0){{\footnotesize $\circledast$}}}
\put(120,130){\makebox(0,0){{\footnotesize $\circledast$}}}
\put(-120,105){\makebox(0,0){\tiny $\bullet$}}
\put(0,105){\makebox(0,0){{\footnotesize $\circledast$}}}
\put(120,105){\makebox(0,0){{\footnotesize $\circledast$}}}
\put(-60,80){\makebox(0,0){{\tiny $\bullet$}}}
\put(0,80){\makebox(0,0){{\tiny $\bullet$}}}
\put(60,80){\makebox(0,0){{\footnotesize $\circledast$}}}
\put(120,80){\makebox(0,0){{\footnotesize $\circledast$}}}

\put(-30,55){\makebox(0,0){{\tiny $\bullet$}}}
\put(30,55){\makebox(0,0){{\tiny $\bullet$}}}
\put(60,55){\makebox(0,0){{\footnotesize $\circledast$}}}
\put(90,55){\makebox(0,0){{\footnotesize $\circledast$}}}
\put(120,55){\makebox(0,0){{\tiny $\bullet$}}}
\put(55.5,-20){\makebox(0,0){{\tiny $\bullet$}}}
\put(64.5,-20){\makebox(0,0){{\tiny $\bullet$}}}
\put(67.5,-20){\makebox(0,0){{\footnotesize $\circledast$}}}
\put(71,-20){\makebox(0,0){{\footnotesize $\circledast$}}}
\put(75,-20){\makebox(0,0){{\tiny $\bullet$}}}
\put(79.5,-20){\makebox(0,0){{\tiny $\bullet$}}}

\put(105,37){\makebox(0,0){{\small $\frac{4}{5}$}}}
\put(95,30){\makebox(0,0){{\small $\frac{3}{4}$}}}
\put(79,30){\makebox(0,0){{$\mathbf{\frac{5}{7}}$}}}
\put(56,30){\makebox(0,0){{$\mathbf{\frac{2}{3}}$}}}
\put(45,37){\makebox(0,0){{\small $\frac{5}{8}$}}}

\put(105,30){\makebox(0,0){{\tiny $\bullet$}}}
\put(90,30){\makebox(0,0){{\tiny $\bullet$}}}
\put(75,30){\makebox(0,0){{\footnotesize $\circledast$}}}
\put(60,30){\makebox(0,0){{\footnotesize $\circledast$}}}
\put(45,30){\makebox(0,0){{\tiny $\bullet$}}}

\put(84,-1){\makebox(0,0){{\small $\frac{8}{11}$}}}
\put(77,12){\makebox(0,0){{$\mathbf{\frac{5}{7}}$}}}
\put(76,-14){\makebox(0,0){{\small $\frac{5}{7}$}}}
\put(67.5,12){\makebox(0,0){{$\mathbf{\frac{7}{10}}$}}}
\put(58,11){\makebox(0,0){{\small $\frac{2}{3}$}}}
\put(50.5,-1){\makebox(0,0){{\small $\frac{7}{11}$}}}

\put(55.5,-27){\makebox(0,0){{\small $\frac{9}{14}$}}}
\put(83.5,-22){\makebox(0,0){{\small $\frac{13}{18}$}}}

\put(65,-27){\makebox(0,0){{$\mathbf{\frac{7}{10}}$}}}
\put(73,-27){\makebox(0,0){{$\mathbf{\frac{12}{17}}$}}}

\put(82.5,5){\makebox(0,0){{\tiny $\bullet$}}}
\put(75,5){\makebox(0,0){{\footnotesize $\circledast$}}}
\put(67.5,5){\makebox(0,0){{\footnotesize $\circledast$}}}
\put(60,5){\makebox(0,0){{\tiny $\bullet$}}}
\put(52.5,5){\makebox(0,0){{\tiny $\bullet$}}}

\end{picture}
\end{center}
\caption{The diagrams $\GG_\theta=D(\A/I_\theta)=R_2\circ L_2
\circ R_1 \circ L_1 \circ \cdots$ (lighter) and
$\GG\setminus\GG_\theta=D(I_\theta)$ (darker) when
$\theta=[1,2,2,1,1,\ldots]$} \label{Figure7}
\end{figure}

If $j_n=j_n(\theta)$ is the unique index for which
$r(n,j_n)<\theta<r(n,j_n+1)$ (see Figure \ref{Figure7}), then
\begin{equation*}
I_\theta \cap \AA_n =\mbox{\small $\displaystyle
\bigoplus\limits_{\substack{0\leq k\leq 2^n \\
k\neq j_n,j_{n+1}}}$} \M_{q(n,k)} .
\end{equation*}
The vertices of $D(\AA/I_\theta)$ are explicitly  related to the
continued fraction decomposition of $\theta$. For each $r\in
\Q_{(0,1)}$, denote $\operatorname{ht}(r)=\min\{ n:\exists k,\
r(n,k)=r\}$. Let $\frac{p_n}{q_n}$ be the continued fraction
approximations of $\theta$, and $h_n=\operatorname{ht} (
\frac{p_n}{q_n})$. With this notation, the labels of the two
vertices at floor $m$ in $\GG_\theta$ are $\frac{p_n}{q_n}$ and
$\frac{p_{n-1}+(m-h_n)p_n}{q_{n-1}+(m-h_n)q_n}$ whenever $h_n\leq
m<h_{n+1}$.

(ii)  For each $\theta=\frac{p}{q}\in \Q_{(0,1)}$ in lowest terms,
consider the Bratteli subdiagram $\GG_\theta$ of $\GG$ defined by
all vertices $(n,j)$ with $r(n,j)=\theta$ and $(m,i)$ with
$(m,i)\Downarrow (n,j)$. The AF algebra associated to $\GG_\theta$
is clearly isomorphic to $\M_q$. Again, the complement
$\GG\setminus\GG_\theta$ is seen to be a directed and hereditary
Bratteli diagram. Therefore there is a primitive ideal $I_\theta$
in $\AA$ such that $D(I_\theta)=\GG\setminus\GG_\theta$ and
$\AA/I_\theta \simeq \M_q$.

Let $n_0-1=n_0(\theta)-1$ be the largest $n\in \N$ for which there
exists $j=j_n(\theta)$ such that $r(n,j)<\theta<r(n,j+1)$. For $n<
n_0$ define $j_n$ as above. By the choice of $n_0$ and the
properties of the Pascal triangle with repetition, for every
$n\geq n_0$ there is $j_n=j_n(\theta)$ with $r(n,j_n)=\theta$. The
ideal $I_\theta$ is generated by the direct summands
$\M_{q(n_0,j_{n_0}-1)}$, $\M_{q(n_0,j_{n_0}+1)}$ and
$\M_{q(n,c_n)}$, $n<n_0$, that is

\begin{equation*}
I_\theta \cap \AA_n =\begin{cases} \hspace{3pt}
\mbox{\small $\displaystyle \bigoplus\limits_{\substack{0\leq k\leq 2^n \\
k\neq j_n,j_{n+1}}}$} \M_{q(n,k)} & \mbox{\rm if $n<n_0$,} \\
\mbox{\small $\displaystyle \bigoplus\limits_{\substack{0\leq k\leq 2^n \\
k\neq j_n}}$} \M_{q(n,k)} & \mbox{\rm if $n\geq n_0$}.
\end{cases}
\end{equation*}

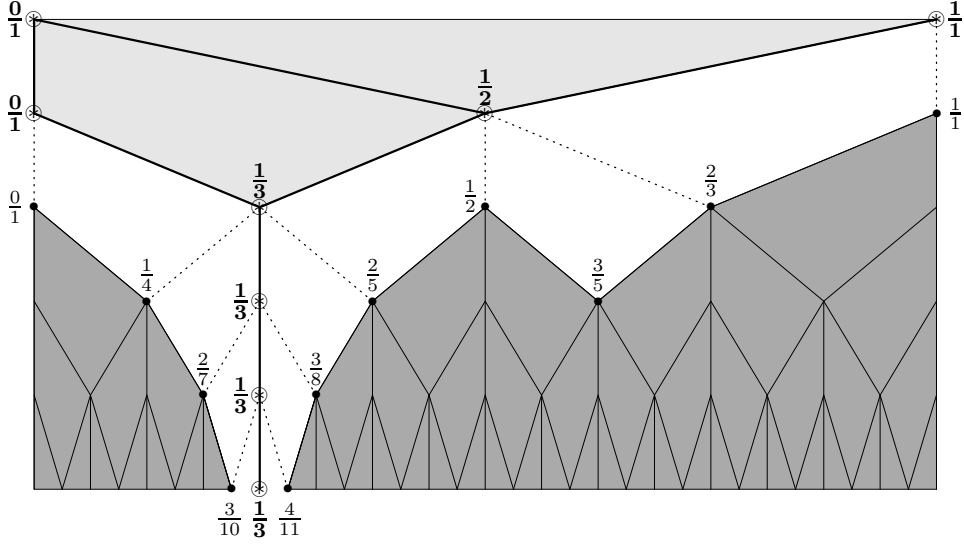
\begin{figure}[htb]
\begin{center}
\unitlength 0.5mm
\begin{picture}(0,135)(5,0)
%\thinlines
\linethickness{10mm} \texture{cc 0}
\shade\path(-120,5)(-120,80)(-90,55)(-75,30)(-67.5,5)(-120,5)
\shade\path(-52.5,5)(-45,30)(-30,55)(0,80)(30,55)(60,80)(120,105)(120,5)(-52.5,5)

\texture{c 0000}
\shade\path(-120,130)(120,130)(0,105)(-60,80)(-120,105)(-120,130)

\thicklines \drawline(-60,5)(-60,80)(-120,105)(-120,130)
\drawline(120,130)(0,105)(-60,80) \drawline(-120,130)(0,105)

\thinlines \dottedline{2}(120,130)(120,105)
\dottedline{2}(60,80)(0,105)(0,80)
\dottedline{2}(-120,105)(-120,80)
\dottedline{2}(-90,55)(-60,80)(-30,55)
\dottedline{2}(-75,30)(-60,55)(-45,30)
\dottedline{2}(-67.5,5)(-60,30)(-52.5,5)

\drawline(120,130)(0,105) \drawline(-105,30)(-105,5)
\drawline(-90,55)(-90,5) \drawline(-75,30)(-75,5)
\drawline(-45,30)(-45,5) \drawline(-30,55)(-30,5)
\drawline(-15,30)(-15,5) \drawline(0,80)(0,5)
\drawline(15,30)(15,5) \drawline(30,55)(30,5)
\drawline(45,30)(45,5) \drawline(60,80)(60,5)
\drawline(75,30)(75,5) \drawline(90,55)(90,5)
\drawline(105,30)(105,5)
\drawline(-120,30)(-112.5,5)(-105,30)(-97.5,5)(-90,30)(-82.5,5)(-75,30)(-67.5,5)
\drawline(-52.5,5)(-45,30)(-37.5,5)(-30,30)(-22.5,5)(-15,30)(-7.5,5)(0,30)
\drawline(120,30)(112.5,5)(105,30)(97.5,5)(90,30)(82.5,5)(75,30)(67.5,5)(60,30)
\drawline(60,30)(52.5,5)(45,30)(37.5,5)(30,30)(22.5,5)(15,30)(7.5,5)(0,30)
\drawline(-120,55)(-105,30)(-90,55)
\drawline(-30,55)(-15,30)(0,55)(15,30)(30,55)(45,30)(60,55)(75,30)(90,55)(105,30)(120,55)
\drawline(60,80)(90,55)(120,80)

\drawline(-120,5)(-120,80)(-90,55)(-75,30)(-67.5,5)
\drawline(120,5)(120,105)(60,80)(30,55)(0,80)(-30,55)(-45,30)(-52.5,5)

\put(-125,130){\makebox(0,0){{$\mathbf{\frac{0}{1}}$}}}
\put(0,112){\makebox(0,0){{$\mathbf{\frac{1}{2}}$}}}
\put(125,130){\makebox(0,0){{$\mathbf{\frac{1}{1}}$}}}
\put(-125,105){\makebox(0,0){{$\mathbf{\frac{0}{1}}$}}}
\put(125,105){\makebox(0,0){{\small $\frac{1}{1}$}}}
\put(-125,80){\makebox(0,0){{\small $\frac{0}{1}$}}}
\put(-60,87){\makebox(0,0){{$\mathbf{\frac{1}{3}}$}}}
\put(-4,83){\makebox(0,0){{\small $\frac{1}{2}$}}}
\put(60,87){\makebox(0,0){{\small $\frac{2}{3}$}}}
%\put(125,80){\makebox(0,0){{\small $\frac{1}{1}$}}}
%\put(-125,80){\makebox(0,0){{\footnotesize $\frac{0}{1}$}}}

\put(-90,62){\makebox(0,0){{\small $\frac{1}{4}$}}}
\put(-65,55){\makebox(0,0){{$\mathbf{\frac{1}{3}}$}}}
\put(-30,62){\makebox(0,0){{\small $\frac{2}{5}$}}}
%\put(-5,55){\makebox(0,0){{\small $\frac{1}{2}$}}}
\put(30,62){\makebox(0,0){{\small $\frac{3}{5}$}}}

\put(-120,130){\makebox(0,0){{\footnotesize $\circledast$}}}
\put(120,130){\makebox(0,0){{\footnotesize $\circledast$}}}
\put(-120,105){\makebox(0,0){{\footnotesize $\circledast$}}}
\put(0,105){\makebox(0,0){{\footnotesize $\circledast$}}}
\put(120,105){\makebox(0,0){{\tiny $\bullet$}}}
\put(-120,80){\makebox(0,0){{\tiny $\bullet$}}}
\put(-60,80){\makebox(0,0){{\footnotesize $\circledast$}}}
\put(0,80){\makebox(0,0){{\tiny $\bullet$}}}
\put(60,80){\makebox(0,0){{\tiny $\bullet$}}}

\put(-90,55){\makebox(0,0){{\tiny $\bullet$}}}
\put(-60,55){\makebox(0,0){{\footnotesize $\circledast$}}}
\put(-30,55){\makebox(0,0){{\tiny $\bullet$}}}
\put(30,55){\makebox(0,0){{\tiny $\bullet$}}}

\put(-75,37){\makebox(0,0){{\small $\frac{2}{7}$}}}
\put(-65,30){\makebox(0,0){{$\mathbf{\frac{1}{3}}$}}}
\put(-45,37){\makebox(0,0){{\small $\frac{3}{8}$}}}

\put(-75,30){\makebox(0,0){{\tiny $\bullet$}}}
\put(-60,30){\makebox(0,0){{\footnotesize $\circledast$}}}
\put(-45,30){\makebox(0,0){{\tiny $\bullet$}}}

\put(-68,-3){\makebox(0,0){{\small $\frac{3}{10}$}}}
\put(-60,-3){\makebox(0,0){{$\mathbf{\frac{1}{3}}$}}}
\put(-52,-3){\makebox(0,0){{\small $\frac{4}{11}$}}}

\put(-67.5,5){\makebox(0,0){{\tiny $\bullet$}}}
\put(-60,5){\makebox(0,0){{\footnotesize $\circledast$}}}
\put(-52.5,5){\makebox(0,0){{\tiny $\bullet$}}}

\end{picture}
\end{center}
\caption{The diagrams $D(I_{\frac{1}{3}})$ (darker) and
$D(\AA/I_{\frac{1}{3}})$ (lighter)} \label{Figure8}
\end{figure}

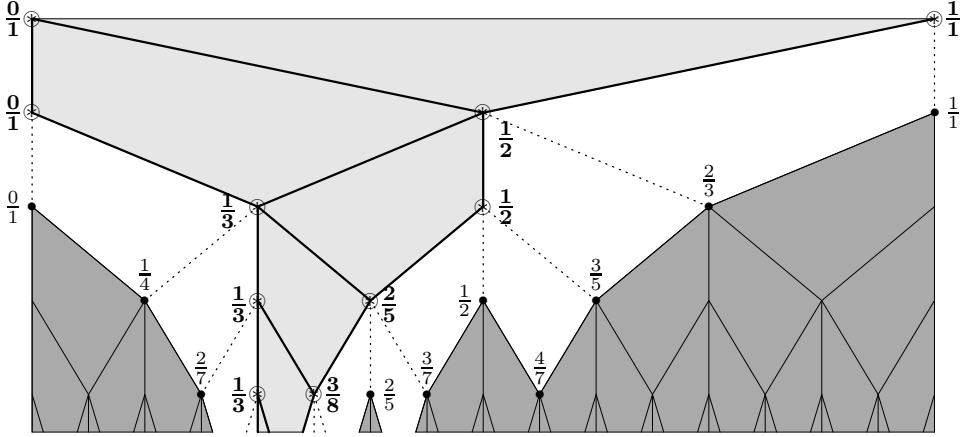
\begin{figure}[htb]
\begin{center}
\unitlength 0.5mm
\begin{picture}(0,120)(5,22)
%\thinlines
\linethickness{10mm} \texture{cc 0}
\shade\path(-120,20)(-120,80)(-90,55)(-75,30)(-72,20)(-120,20)
\shade\path(120,105)(120,20)(-18,20)(-15,30)(0,55)(15,30)(30,55)(60,80)(120,105)
\shade\path(-30,30)(-27,20)(-33,20)(-30,30)

\texture{c 0000}
\shade\path(-120,130)(120,130)(0,105)(0,80)(-30,55)(-45,30)(-48,20)(-60,20)(-60,80)(-120,105)(-120,130)

\thicklines \drawline(-120,130)(0,105)(-60,80)(-30,55)
\drawline(-60,55)(-45,30) \drawline(-60,30)(-57,20)
\drawline(-48,20)(-45,30)(-30,55)(0,80)(0,105)(120,130)
\drawline(-120,130)(-120,105)(-60,80)(-60,20)

\thinlines \dottedline{2}(120,130)(120,105)
\dottedline{2}(60,80)(0,105) \dottedline{2}(0,55)(0,80)(30,55)
\dottedline{2}(-120,105)(-120,80) \dottedline{2}(-90,55)(-60,80)
\dottedline{2}(-75,30)(-60,55) \dottedline{2}(-60,30)(-63,20)
\dottedline{2}(-45,20)(-45,30)(-42,20)
\dottedline{2}(-30,30)(-30,55)(-15,30)

\drawline(-120,20)(-120,80)(-90,55)(-75,30)(-72,20)
\drawline(120,20)(120,105)(60,80)(30,55)(15,30)(0,55)(-15,30)(-18,20)
\drawline(-27,20)(-30,30)(-33,20) \drawline(-105,30)(-105,20)
\drawline(-108,20)(-105,30)(-102,20) \drawline(-90,55)(-90,20)
\drawline(-93,20)(-90,30)(-87,20)
\drawline(-78,20)(-75,30)(-75,20)
\drawline(-120,55)(-105,30)(-90,55) \drawline(-30,30)(-30,20)
\drawline(-15,20)(-15,30)(-12,20) \drawline(0,55)(0,20)
\drawline(-3,20)(0,30)(3,20) \drawline(15,30)(15,20)
\drawline(12,20)(15,30)(18,20) \drawline(30,55)(30,20)
\drawline(27,20)(30,30)(33,20) \drawline(45,30)(45,20)
\drawline(42,20)(45,30)(48,20) \drawline(75,30)(75,20)
\drawline(72,20)(75,30)(78,20) \drawline(105,30)(105,20)
\drawline(102,20)(105,30)(108,20) \drawline(60,80)(60,20)
\drawline(90,55)(90,20) \drawline(120,80)(90,55)(60,80)
\drawline(30,55)(45,30)(60,55)(75,30)(90,55)(105,30)(120,55)
\drawline(57,20)(60,30)(63,20) \drawline(87,20)(90,30)(93,20)
\drawline(120,30)(117,20) \drawline(-120,30)(-117,20)

\put(-125,130){\makebox(0,0){{$\mathbf{\frac{0}{1}}$}}}
\put(6,98){\makebox(0,0){{$\mathbf{\frac{1}{2}}$}}}
\put(125,130){\makebox(0,0){{$\mathbf{\frac{1}{1}}$}}}
\put(-125,105){\makebox(0,0){{$\mathbf{\frac{0}{1}}$}}}
\put(125,105){\makebox(0,0){{\small $\frac{1}{1}$}}}
\put(-125,80){\makebox(0,0){{\small $\frac{0}{1}$}}}
\put(-68,79){\makebox(0,0){{$\mathbf{\frac{1}{3}}$}}}
\put(6,80){\makebox(0,0){{$\mathbf{\frac{1}{2}}$}}}
\put(60,87){\makebox(0,0){{\small $\frac{2}{3}$}}}
%\put(125,80){\makebox(0,0){{\small $\frac{1}{1}$}}}

%\put(-125,55){\makebox(0,0){{\small $\frac{0}{1}$}}}
\put(-90,62){\makebox(0,0){{\small $\frac{1}{4}$}}}
\put(-65,54){\makebox(0,0){{$\mathbf{\frac{1}{3}}$}}}
\put(-25,54){\makebox(0,0){{$\mathbf{\frac{2}{5}}$}}}
\put(-5,55){\makebox(0,0){{\small $\frac{1}{2}$}}}
\put(30,62){\makebox(0,0){{\small $\frac{3}{5}$}}}

\put(-120,130){\makebox(0,0){{\footnotesize $\circledast$}}}
\put(120,130){\makebox(0,0){{\footnotesize $\circledast$}}}
\put(-120,105){\makebox(0,0){{\footnotesize $\circledast$}}}
\put(0,105){\makebox(0,0){{\footnotesize $\circledast$}}}
\put(120,105){\makebox(0,0){{\tiny $\bullet$}}}
\put(-120,80){\makebox(0,0){{\tiny $\bullet$}}}
\put(-60,80){\makebox(0,0){{\footnotesize $\circledast$}}}
\put(0,80){\makebox(0,0){{\footnotesize $\circledast$}}}
\put(60,80){\makebox(0,0){{\tiny $\bullet$}}}

\put(-90,55){\makebox(0,0){{\tiny $\bullet$}}}
\put(-60,55){\makebox(0,0){{\footnotesize $\circledast$}}}
\put(-30,55){\makebox(0,0){{\footnotesize $\circledast$}}}
\put(0,55){\makebox(0,0){{\tiny $\bullet$}}}
\put(30,55){\makebox(0,0){{\tiny $\bullet$}}}

\put(-75,37){\makebox(0,0){{\small $\frac{2}{7}$}}}
\put(-65,30){\makebox(0,0){{$\mathbf{\frac{1}{3}}$}}}
\put(-40,30){\makebox(0,0){{$\mathbf{\frac{3}{8}}$}}}
\put(-25,30){\makebox(0,0){{\small $\frac{2}{5}$}}}
\put(-15,37){\makebox(0,0){{\small $\frac{3}{7}$}}}
\put(15,37){\makebox(0,0){{\small $\frac{4}{7}$}}}

\put(-75,30){\makebox(0,0){{\tiny $\bullet$}}}
\put(-60,30){\makebox(0,0){{\footnotesize $\circledast$}}}
\put(-45,30){\makebox(0,0){{\footnotesize $\circledast$}}}
\put(-30,30){\makebox(0,0){{\tiny $\bullet$}}}
\put(-15,30){\makebox(0,0){{\tiny $\bullet$}}}
\put(15,30){\makebox(0,0){{\tiny $\bullet$}}}

\end{picture}
\end{center}
\caption{The diagrams $D(I^+_{\frac{1}{3}})$ (darker) and
$D(\AA/I^+_{\frac{1}{3}})$ (lighter)} \label{Figure9}
\end{figure}

\begin{figure}[htb]
\begin{center}
\unitlength 0.5mm
\begin{picture}(0,135)(5,0)
%\thinlines
\linethickness{10mm} \texture{cc 0}
\shade\path(-43,-10)(-45,0)(-47,-10)(-43,-10)
\shade\path(-120,-10)(-120,75)(-90,50)(-75,25)(-67.5,0)(-60,25)(-52.5,0)(-50.5,-10)(-120,-10)
\shade\path(120,-10)(120,100)(60,75)(30,50)(15,25)(0,50)(-15,25)(-22.5,0)(-24.5,-10)(120,-10)

\texture{c 0000}
\shade\path(120,125)(-120,125)(-120,100)(-60,75)(-60,50)(-45,25)(-37.5,0)(-35.5,-10)(-30,-10)(-30,50)(0,75)(0,100)(120,125)

\thicklines \drawline(-120,125)(0,100)(-60,75)(-30,50)(-45,25)
\drawline(-30,25)(-37.5,0) \drawline(-30,0)(-32,-10)
\drawline(-30,-10)(-30,50)(0,75)(0,100)(120,125)
\drawline(-120,125)(-120,100)(-60,75)(-60,50)(-45,25)(-37.5,0)(-35.5,-10)

\thinlines \dottedline{2}(120,125)(120,100)
\dottedline{2}(-120,75)(-120,100) \dottedline{2}(-90,50)(-60,75)
\dottedline{2}(-75,25)(-60,50)(-60,25)
\dottedline{2}(-52.5,0)(-45,25)(-45,0)
\dottedline{2}(-39.5,-10)(-37.5,0)(-37.5,-10)
\dottedline{2}(120,125)(120,100) \dottedline{2}(60,75)(0,100)
\dottedline{2}(30,50)(0,75)(0,50) \dottedline{2}(-30,50)(-15,25)
\dottedline{2}(-30,25)(-22.5,0) \dottedline{2}(-30,0)(-28,-10)

%\drawline(-45,0)(-45,-10)
\drawline(-105,25)(-105,0)
%\drawline(-112.5,0)(-112.5,-10) \drawline(-120,0)(-118,-10)
%\drawline(-114.5,-10)(-112.5,5)(-110.5,-5)
\drawline(-120,50)(-105,25)(-90,50)(-90,0)
\drawline(-120,25)(-112.5,0)(-105,25)(-97.5,0)(-90,25)(-82.5,0)(-75,25)
%\drawline(-103,-5)(-105,5)(-107,-5)
%\drawline(-99.5,-5)(-97.5,5)(-95.5,-5)
%\drawline(-92,-5)(-90,5)(-88,-5) \drawline(-97.5,5)(-97.5,-5)
%\drawline(-82.5,5)(-82.5,-5)
\drawline(-75,25)(-75,0) %\drawline(-84.5,-5)(-82.5,5)(-80.5,-5)
%\drawline(-77,-5)(-75,5)(-73,-5) \drawline(-67.5,5)(-67.5,-5)
%\drawline(-69.5,-5)(-67.5,5)(-65.5,-5)
\drawline(-60,25)(-60,0) %\drawline(-62,-5)(-60,5)(-58,-5)
%\drawline(-54.5,-5)(-52.5,5)(-52.5,-5) \drawline(120,5)(118,-5)
%\drawline(112.5,5)(112.5,-5)
%\drawline(-114.5,-5)(-112.5,5)(-110.5,-5)
%\drawline(-22.5,-5)(-22.5,5)(-20.5,-5) \drawline(-7.5,5)(-7.5,-5)
%\drawline(-9.5,-5)(-7.5,5)(-5.5,-5)
\drawline(-15,25)(-15,0) %\drawline(-17,-5)(-15,5)(-13,-5)
\drawline(0,50)(0,0) %\drawline(-2,-5)(0,5)(2,-5)
%\drawline(7.5,5)(7.5,-5) \drawline(5.5,-5)(7.5,5)(9.5,-5)
\drawline(15,25)(15,0) %\drawline(13,-5)(15,5)(17,-5)
%\drawline(22.5,5)(22.5,-5) \drawline(20.5,-5)(22.5,5)(24.5,-5)
\drawline(30,50)(30,0) %\drawline(28,-5)(30,5)(32,-5)
%\drawline(37.5,5)(37.5,-5) \drawline(35.5,-5)(37.5,5)(39.5,-5)
\drawline(45,25)(45,0) %\drawline(43,-5)(45,5)(47,-5)
%\drawline(52.5,-5)(52.5,5) \drawline(50.5,-5)(52.5,5)(54.5,-5)
\drawline(60,75)(60,0) %\drawline(58,-5)(60,5)(62,-5)
%\drawline(67.5,5)(67.5,-5) \drawline(65.5,-5)(67.5,5)(69.5,-5)
\drawline(75,25)(75,0) %\drawline(73,-5)(75,5)(77,-5)
%\drawline(82.5,5)(82.5,-5) \drawline(80.5,-5)(82.5,5)(84.5,-5)
\drawline(90,50)(90,0) %\drawline(88,-5)(90,5)(92,-5)
%\drawline(97.5,5)(97.5,-5) \drawline(95.5,-5)(97.5,5)(99.5,-5)
\drawline(105,25)(105,0) %\drawline(103,-5)(105,5)(107,-5)
%\drawline(110.5,-5)(112.5,5)(114.5,-5)
\drawline(-15,25)(-7.5,0)(0,25)(7.5,0)(15,25)(22.5,0)(30,25)(37.5,0)(45,25)(52.5,0)(60,25)(67.5,0)(75,25)(82.5,0)(90,25)(97.5,0)(105,25)(112.5,0)(120,25)
\drawline(30,50)(45,25)(60,50)(75,25)(90,50)(105,25)(120,50)
\drawline(60,75)(90,50)(120,75)

%\drawline(-43,-5)(-45,5)(-47,-5)
\drawline(-120,0)(-120,75)(-90,50)(-75,25)(-67.5,0)(-60,25)(-52.5,0)
\drawline(120,0)(120,100)(60,75)(30,50)(15,25)(0,50)(-15,25)(-22.5,0)

\put(-125,125){\makebox(0,0){{$\mathbf{\frac{0}{1}}$}}}
\put(6,96){\makebox(0,0){{$\mathbf{\frac{1}{2}}$}}}
\put(125,125){\makebox(0,0){{$\mathbf{\frac{1}{1}}$}}}
\put(-125,100){\makebox(0,0){{$\mathbf{\frac{0}{1}}$}}}
\put(125,100){\makebox(0,0){{\small $\frac{1}{1}$}}}
\put(-125,75){\makebox(0,0){{\small $\frac{0}{1}$}}}
\put(-68,74){\makebox(0,0){{$\mathbf{\frac{1}{3}}$}}}
\put(6,75){\makebox(0,0){{$\mathbf{\frac{1}{2}}$}}}
\put(60,82){\makebox(0,0){{\small $\frac{2}{3}$}}}

\put(-90,57){\makebox(0,0){{\small $\frac{1}{4}$}}}
\put(-65,50){\makebox(0,0){{$\mathbf{\frac{1}{3}}$}}}
\put(-25,49){\makebox(0,0){{$\mathbf{\frac{2}{5}}$}}}
\put(-5,50){\makebox(0,0){{\small $\frac{1}{2}$}}}
\put(30,57){\makebox(0,0){{\small $\frac{3}{5}$}}}

\put(-120,125){\makebox(0,0){{\footnotesize $\circledast$}}}
\put(120,125){\makebox(0,0){{\footnotesize $\circledast$}}}
\put(-120,100){\makebox(0,0){{\footnotesize $\circledast$}}}
\put(0,100){\makebox(0,0){{\footnotesize $\circledast$}}}
\put(120,100){\makebox(0,0){{\tiny $\bullet$}}}
\put(-120,75){\makebox(0,0){{\tiny $\bullet$}}}
\put(-60,75){\makebox(0,0){{\footnotesize $\circledast$}}}
\put(0,75){\makebox(0,0){{\footnotesize $\circledast$}}}
\put(60,75){\makebox(0,0){{\tiny $\bullet$}}}
%\put(120,75){\makebox(0,0){{\tiny $\bullet$}}}

%\put(-120,50){\makebox(0,0){{\tiny $\bullet$}}}
\put(-90,50){\makebox(0,0){{\tiny $\bullet$}}}
\put(-60,50){\makebox(0,0){{\footnotesize $\circledast$}}}
\put(-30,50){\makebox(0,0){{\footnotesize $\circledast$}}}
\put(0,50){\makebox(0,0){{\tiny $\bullet$}}}
\put(30,50){\makebox(0,0){{\tiny $\bullet$}}}

\put(-75,32){\makebox(0,0){{\small $\frac{2}{7}$}}}
\put(-65,25){\makebox(0,0){{\small $\frac{1}{3}$}}}
\put(-50,25){\makebox(0,0){{$\mathbf{\frac{3}{8}}$}}}
\put(-26,25){\makebox(0,0){{\small $\frac{2}{5}$}}}
\put(-15,32){\makebox(0,0){{\small $\frac{3}{7}$}}}
\put(15,32){\makebox(0,0){{\small $\frac{4}{7}$}}}

\put(-75,25){\makebox(0,0){{\tiny $\bullet$}}}
\put(-60,25){\makebox(0,0){{\tiny $\bullet$}}}
\put(-45,25){\makebox(0,0){{\footnotesize $\circledast$}}}
\put(-30,25){\makebox(0,0){{\footnotesize $\circledast$}}}
\put(-15,25){\makebox(0,0){{\tiny $\bullet$}}}
\put(15,25){\makebox(0,0){{\tiny $\bullet$}}}

\put(-67.5,13){\makebox(0,0){{\footnotesize $\frac{3}{10}$}}}
\put(-51.5,12){\makebox(0,0){{\footnotesize $\frac{4}{11}$}}}
\put(-42,0){\makebox(0,0){{\small $\frac{3}{8}$}}}
\put(-36,8){\makebox(0,0){{\small $\mathbf{\frac{5}{13}}$}}}
\put(-27,0){\makebox(0,0){{\footnotesize $\frac{2}{5}$}}}
\put(-23.5,12){\makebox(0,0){\footnotesize $\frac{5}{12}$}}

\put(-67.5,0){\makebox(0,0){{\tiny $\bullet$}}}
\put(-52.5,0){\makebox(0,0){{\tiny $\bullet$}}}
\put(-45,0){\makebox(0,0){{\tiny $\bullet$}}}
\put(-30,0){\makebox(0,0){{\footnotesize $\circledast$}}}
\put(-22.5,0){\makebox(0,0){{\tiny $\bullet$}}}
\put(-37.5,0){\makebox(0,0){{\footnotesize $\circledast$}}}

\end{picture}
\end{center}
\caption{The diagrams $D(I_{\frac{2}{5}}^-)$ (darker) and
$D(\AA/I^-_{\frac{2}{5}})$ (lighter)} \label{Figure10}
\end{figure}
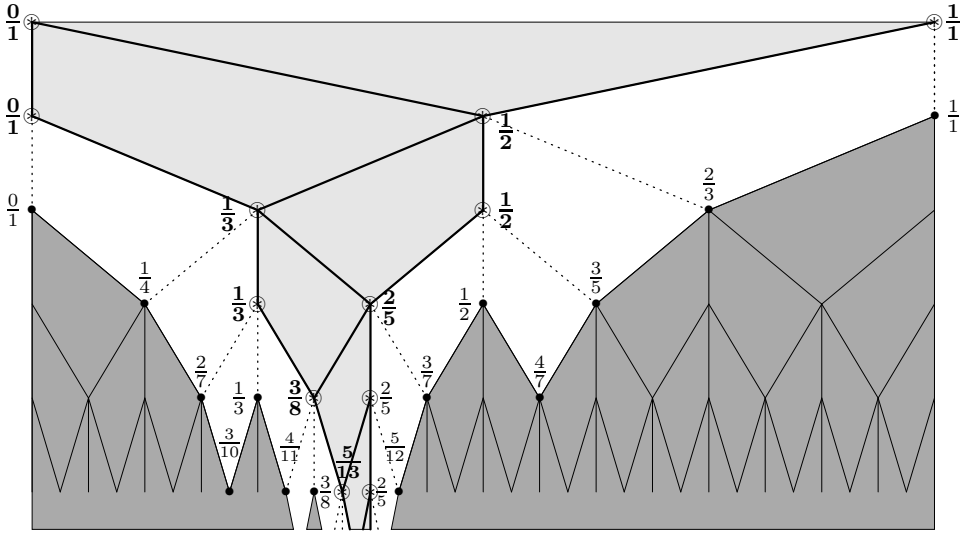

The ideals $I_\theta^\pm$ defined by (see also Figures
\ref{Figure9} and \ref{Figure10})
\begin{equation*}
I_\theta^+ \cap \AA_n =\mbox{\small $\displaystyle
\bigoplus_{\substack{0\leq k\leq 2^n \\ k\neq j_n,j_{n+1}}}$}
\M_{q(n,k)},
\end{equation*}
and respectively by
\begin{equation*}
I_\theta^- \cap \AA_n =\begin{cases} \hspace{3pt}
\mbox{\small $\displaystyle \bigoplus\limits_{\substack{0\leq k\leq 2^n \\
k\neq j_n,j_{n+1}}}$} \M_{q(n,k)} & \mbox{\rm if $n<n_0$,} \\
\mbox{\small $\displaystyle \bigoplus\limits_{\substack{0\leq k\leq 2^n \\
k\neq j_{n-1},j_n}}$} \M_{q(n,k)} & \mbox{\rm if $n\geq n_0$},
\end{cases}
\end{equation*}
are primitive and we clearly have $\AA/I_\theta^- \cong
\AA_{(q,\overline{p})}$ and $\AA/I^+_\theta \cong
\AA_{(q,q-\overline{p})}$.

(iii) is now obvious.
\end{proof}

\begin{remark}\label{R3}
{\em A joint (and important) feature of all cases above is that}
\begin{equation*}
(n,j)\notin D(I_\theta)=\GG\setminus\GG_\theta \ \Longrightarrow \
r(n,j-1)<\theta <r(n,j+1).
\end{equation*}
\end{remark}

\begin{remark}\label{R4}
{\em In $GL_2(\Z)$ consider the matrices
\begin{equation*}
A=\left[ \begin{matrix} 1 & 0 \\ 1 & 1 \end{matrix} \right],\quad
B=\left[ \begin{matrix} 1 & 1 \\ 0 & 1 \end{matrix}\right],\quad
J=\left[\begin{matrix} 0 & 1 \\ 1 & 0 \end{matrix} \right] ,\quad
M(a)=\left[ \begin{matrix} a & 1 \\ 1 & 0 \end{matrix} \right] .
\end{equation*}
The identification between $L_a \circ R_b$ and $C_a \circ C_b$
reflects the matrix equality \begin{equation*} B^a A^b =M(a)M(b),
\end{equation*}
whereas the identification between $R_a \circ R_b$ and $C_a \circ
C_b$ reflects the matrix equality
\begin{equation*}
A^a B^b=JM(a)M(b)J.
\end{equation*} }
\end{remark}

A combinatorial analysis based on Bratteli's correspondence
between primitive ideals and subdiagrams of $\GG$ shows that these
are actually the only primitive ideals of $\AA$.

\begin{proposition}\label{P5}
$\Prim \AA=\{ I_\theta:\theta \in \I\} \cup \{
I_\theta,I_\theta^\pm :\theta \in \Q_{(0,1)}\} \cup \{
I_0,I_0^+,I_1,I_1^-\}$.
\end{proposition}

\begin{figure}[htb]
\begin{center}
\unitlength 0.3mm
\begin{picture}(70,80)(-10,-38)
\dottedline{3}(-160,10)(-160,40)(-130,10)(-100,40)(-70,10)(-40,40)(-40,10)
\path(-100,10)(-100,40) \put(-100,10){\makebox(0,0){{\small
$\circledast$}}} \put(-100,40){\makebox(0,0){{\small
$\circledast$}}} \put(-160,10){\makebox(0,0){{\tiny $\bullet$}}}
\put(-160,40){\makebox(0,0){{\tiny $\bullet$}}}
\put(-130,10){\makebox(0,0){{\tiny $\bullet$}}}
\put(-70,10){\makebox(0,0){{\tiny $\bullet$}}}
\put(-40,40){\makebox(0,0){{\tiny $\bullet$}}}
\put(-40,10){\makebox(0,0){{\tiny $\bullet$}}}

\dottedline{3}(20,40)(50,10)(80,40)
\path(80,10)(80,40)(110,10)(140,40)
\dottedline{3}(140,10)(140,40)(170,10)(200,40)
\put(140,40){\makebox(0,0){{\small $\circledast$}}}
\put(80,40){\makebox(0,0){{\small $\circledast$}}}
\put(80,10){\makebox(0,0){{\small $\circledast$}}}
\put(110,10){\makebox(0,0){{\small $\circledast$}}}
\put(50,10){\makebox(0,0){{\tiny $\bullet$}}}
\put(20,40){\makebox(0,0){{\tiny $\bullet$}}}
\put(140,10){\makebox(0,0){{\tiny $\bullet$}}}
\put(170,10){\makebox(0,0){{\tiny $\bullet$}}}
\put(200,40){\makebox(0,0){{\tiny $\bullet$}}}

\dottedline{3}(-170,-10)(-140,-40)(-110,-10)(-110,-40)
\path(-110,-10)(-80,-40)(-50,-10)(-50,-40)
\dottedline{3}(-50,-10)(-20,-40)(10,-10)
\put(-50,-40){\makebox(0,0){{\small $\circledast$}}}
\put(-50,-10){\makebox(0,0){{\small $\circledast$}}}
\put(-80,-40){\makebox(0,0){{\small $\circledast$}}}
\put(-110,-10){\makebox(0,0){{\small $\circledast$}}}
\put(-110,-40){\makebox(0,0){{\tiny $\bullet$}}}
\put(-140,-40){\makebox(0,0){{\tiny $\bullet$}}}
\put(-170,-10){\makebox(0,0){{\tiny $\bullet$}}}
\put(10,-10){\makebox(0,0){{\tiny $\bullet$}}}
\put(-20,-40){\makebox(0,0){{\tiny $\bullet$}}}

\dottedline{3}(40,-10)(70,-40)(100,-10)(100,-40)
\path(100,-10)(130,-40)(160,-10)
\dottedline{3}(160,-40)(160,-10)(190,-40)(220,-10)
\put(100,-10){\makebox(0,0){{\small $\circledast$}}}
\put(160,-10){\makebox(0,0){{\small $\circledast$}}}
\put(130,-40){\makebox(0,0){{\small $\circledast$}}}
\put(160,-40){\makebox(0,0){{\tiny $\bullet$}}}
\put(190,-40){\makebox(0,0){{\tiny $\bullet$}}}
\put(220,-10){\makebox(0,0){{\tiny $\bullet$}}}
\put(40,-10){\makebox(0,0){{\tiny $\bullet$}}}
\put(70,-40){\makebox(0,0){{\tiny $\bullet$}}}
\put(100,-40){\makebox(0,0){{\tiny $\bullet$}}}

\end{picture}
\end{center}
\caption{The possible links between two consecutive floors in
$D(\AA/I)$} \label{Figure11}
\end{figure}
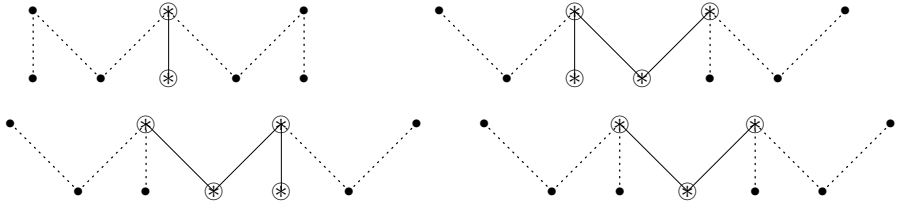

\begin{proof}
Let $I\in \Prim \AA$. Consider the Bratteli diagrams $D=D(I)$ and
$\widetilde{D}=D(\AA/I)=\GG\setminus D$. If there is $n_0$ such
that $(n_0,k)\in D$ for all $0\leq k\leq 2^{n_0}$, then $I=\AA$.
So for each $n$ the set $L_n=\{ k:(n,k)\in \widetilde{D}\}$ is
nonempty. Denote also $L_n^c=\{ 0,1,\ldots,2^n\} \setminus L_n$.

We first notice that $L_n$ should be a set of the form $\{ a_n\}$
or $\{ a_n,a_n+1\}$. If not, there are $k,k^\prime \in L_n$ such
that $k^\prime -k\geq 2$. Since $I$ is a primitive ideal, a vertex
$(p,r)$ in $\GG$ should exist such that $(n,k)\Downarrow (p,r)$
and $(n,k')\Downarrow (p,r)$. Since $k'-k>2$ this is not possible
due to the definition of $\GG$.

To finish the proof it suffices to show that
\begin{equation}\label{2.3}
L_{n+1}=\begin{cases} \{2a_n\} & \mbox{\rm if $L_n=\{ a_n\},$} \\
\{ 2a_n,2a_n+1\},\ \{ 2a_n+1,2a_n+2\}, & \\\quad \mbox{\rm or}\ \
\{ 2a_n+1\} & \mbox{\rm if $L_n=\{ a_n,a_n+1\}$,}
\end{cases}
\end{equation}
that is, all links $(n,j)\downarrow (n+1,j^\prime)$ in
$\widetilde{D}$ are exactly as indicated in Figure \ref{Figure11}.

Indeed, if $L_n=\{ a_n\}$, then $(n,a_n-1),(n,a_n+1)$ are vertices
in the hereditary diagram $D$; thus we also have
$(n+1,2a_n-1),(n+1,2a_n+1)\in D$. Because $D$ is directed,
$(n+1,2a_n)\in D$ would imply $(n,a_n)\in D$, which contradicts
$a_n\in L_n$.

If $L_n=\{ a_n,a_n+1\}$, then $(n,a_n-1),(n,a_n+2)\in D$. Moreover
because $D$ is hereditary the vertices $(n+1,2a_n-1)$ and
$(n+1,2a_n+3)$ also belong to $D$. We now look at the consecutive
vertices $(n+1,2a_n),(n+1,2a_n+1)$, $(n+1,2a_n+2)$. From the first
part they cannot all belong to $\widetilde{D}$. If
$(n+1,2a_n+1)\in D$, and $(n+1,2a_n),(n+1,2a_n+2)\in
\widetilde{D}$, then $L_{n+1}$ has a gap, thus contradicting the
first part. If $(n+1,2a_n),(n+1,2a_n+2)\in D$ it follows, as a
result of the fact that $(n+1,2a_n-1)\in D$ and that $D$ is
directed, that $(n+1,2a_n +1)\in \widetilde{D}$. In a similar way
one cannot have $(n+1,2a_n+1),(n+1,2a_n+2)\in D$. It remains that
only the following cases can occur (see also Figure
\ref{Figure11}):
\begin{itemize}
\item[(i)] $(n+1,2a_n),(n+1,2a_n+1)\in \widetilde{D}$ and
$(n+1,2a_n+2)\in D$, thus $L_{n+1}=\{ 2a_n,2a_n+1\}$. \item[(ii)]
$(n+1,2a_n)\in D$ and $(n+1,2a_n+1),(n+1,2a_n+2)\in
\widetilde{D}$, thus $L_{n+1}=\{ 2a_{n}+1,2a_n+2\}$. \item[(iii)]
$(n+1,2a_n+1)\in \widetilde{D}$ and $(n+1,2a_n),(n+1,2a_n+2)\in
D$, thus $L_{n+1}=\{ 2a_n+1\}$,
\end{itemize}
which concludes the proof of \eqref{2.3}.
\end{proof}

\section{The Jacobson topology on $\Prim \AA$}
We first recall some basic things about the primitive ideal space
of a $C^*$-algebra ${\mathcal A}$ following \cite{Dix} and
\cite{RaWi}. For each set $S\subseteq \Prim {\mathcal A}$,
consider the ideal $k(S):= \cap_{J\in S}\, J$ in ${\mathcal A}$,
called the \emph{kernel} of $S$. For each ideal $I$ consider its
\emph{hull}, $h(I):=\{ P\in \Prim {\mathcal A}: I\subseteq P\}$.
The \emph{closure} of a set $S\subseteq \Prim {\mathcal A}$ is
defined as
\begin{equation*}
\overline{S}:=\{ P\in \Prim {\mathcal A}: k(S) \subseteq P\}.
\end{equation*}
There is a unique topology on $\Prim {\mathcal A}$, called the
\emph{Jacobson} (or \emph{hull-kernel}) \emph{topology} such that
its closed sets are exactly those with $S=\overline{S}$. The open
sets in $\Prim {\mathcal A}$ are then precisely those of the form
\begin{equation*}
\OO_I:=\{ P\in \Prim {\mathcal A} : I\nsubseteq P\}
\end{equation*}
for some ideal $I$ in ${\mathcal A}$. The Jacobson topology is
always $T_0$, i.e. for any two distinct points in $\Prim {\mathcal
A}$ one of them has a neighborhood which does not contain the
other.

Moreover, the correspondence $S\mapsto k(S)$ establishes a
one-to-one correspondence between the closed subsets $S$ of $\Prim
{\mathcal A}$ and the lattice of ideals in ${\mathcal A}$, with
inverse given by $I\mapsto h(I)$. For any ideal $I$ in ${\mathcal
A}$, let $p_I$ denote the quotient map ${\mathcal A}\rightarrow
{\mathcal A}/I$. The mapping $P\mapsto P\cap I$ is a homeomorphism
of the open set $\OO_I$ onto $\Prim I$, whereas $Q\mapsto
p_I^{-1}(Q)$ is a homeomorphism of $\Prim {\mathcal A}/I$ onto the
closed set $h(I)$ of $\Prim {\mathcal A}$. A general study of the
primitive ideal space of AF algebras was pursued in
\cite{Bra3,BraEl,Doo}.

We collect some immediate properties of the primitive ideal space
of $\AA$ in the following

\begin{remark}\label{R6}
{\em (i) For each $\theta \in \I$, $\ \overline{ \{ I_\theta
\}}=\{ I_\theta \}$.}

{\em (ii) For each $\theta \in \Q_{(0,1)}$, $I_\theta \nsubseteq
I_\theta^+$, $I_\theta \nsubseteq I_\theta^-$, and $I_\theta
=I_\theta^+ \cap I_\theta^-$. We also have $I_0\nsubseteq I_0^+$
and $I_1 \nsubseteq I_1^-$. Therefore $\overline{ \{ I_\theta
\}}=\{ I_\theta,I_\theta^+,I_\theta^-\}$ whenever $\theta \in
\Q_{(0,1)}$, $\overline{\{ I_0\}}=\{ I_0,I_0^+\}$ and
$\overline{\{ I_1\}}=\{ I_1,I_1^-\}$, showing in particular that
the Jacobson topology on $\Prim \AA$ is not Hausdorff. In spite of
this we shall see that after removing the ``singular points"
$I_\theta^\pm$ from $\Prim \AA$ we retrieve the usual topology on
$[0,1]$.}
\end{remark}

For each set $E\subseteq [0,1]$, consider the ideal
\begin{equation}\label{3.1}
\II(E):=\mbox{\footnotesize $\displaystyle \bigcap_{\theta \in
E}$} I_\theta ,
\end{equation}
and denote by $\overline{E}$ the usual closure of $E$ in $[0,1]$.

\begin{lemma}\label{L7}
$\II(E)=\II(\overline{E})$ for every set $E\subseteq [0,1]$.
\end{lemma}

\begin{figure}[htb]
\begin{center}
\unitlength 1.2mm
\begin{picture}(0,45)
%\thinlines
\linethickness{10mm} \thinlines
\dottedline{1}(-20,0)(-20,40)(20,40)(20,0)

\thicklines \texture{c 000000}
\shade\path(19.25,0)(17.5,10)(15,20)(10,30)(0,40)(-10,30)(-15,20)(-17.5,10)(-19.25,0)(19.25,0)

\put(10,30){\makebox(0,0){{\tiny $\bullet$}}}
\put(15,20){\makebox(0,0){{\tiny $\bullet$}}}
\put(17.5,10){\makebox(0,0){{\tiny $\bullet$}}}
\put(19.25,0){\makebox(0,0){{\tiny $\bullet$}}}
\put(0,40){\makebox(0,0){{\tiny $\bullet$}}}
\put(-19.25,0){\makebox(0,0){{\tiny $\bullet$}}}
\put(-17.5,10){\makebox(0,0){{\tiny $\bullet$}}}
\put(-15,20){\makebox(0,0){{\tiny $\bullet$}}}
\put(-10,30){\makebox(0,0){{\tiny $\bullet$}}}
\put(20,40){\makebox(0,0){{\tiny $\bullet$}}}
\put(20,30){\makebox(0,0){{\tiny $\bullet$}}}
\put(20,20){\makebox(0,0){{\tiny $\bullet$}}}
\put(20,10){\makebox(0,0){{\tiny $\bullet$}}}
\put(20,0){\makebox(0,0){{\tiny $\bullet$}}}
\put(-20,40){\makebox(0,0){{\tiny $\bullet$}}}
\put(-20,30){\makebox(0,0){{\tiny $\bullet$}}}
\put(-20,20){\makebox(0,0){{\tiny $\bullet$}}}
\put(-20,10){\makebox(0,0){{\tiny $\bullet$}}}
\put(-20,0){\makebox(0,0){{\tiny $\bullet$}}}

\put(0,42){\makebox(0,0){\footnotesize $(n,j)$}}
\put(-20,42){\makebox(0,0){\footnotesize $(n,j-1)$}}
\put(20,42){\makebox(0,0){\footnotesize $(n,j+1)$}}
\put(12,32){\makebox(0,0){\footnotesize $(n+1,2j+1)$}}
\put(-12,32){\makebox(0,0){\footnotesize $(n+1,2j-1)$}}
\put(7,17){\makebox(0,0){\footnotesize $(n+1,4j+3)$}}
\put(-5,20){\makebox(0,0){\footnotesize $(n+1,4j-3)$}}
\put(30,30){\makebox(0,0){\footnotesize $(n+1,2j+2)$}}
\put(30,20){\makebox(0,0){\footnotesize $(n+1,4j+4)$}}
\put(-30,30){\makebox(0,0){\footnotesize $(n+1,2j-2)$}}
\put(-30,20){\makebox(0,0){\footnotesize $(n+1,4j-4)$}}

\end{picture}
\end{center}
\caption{The ideal generated by $(n,j)$} \label{Figure12}
\end{figure}
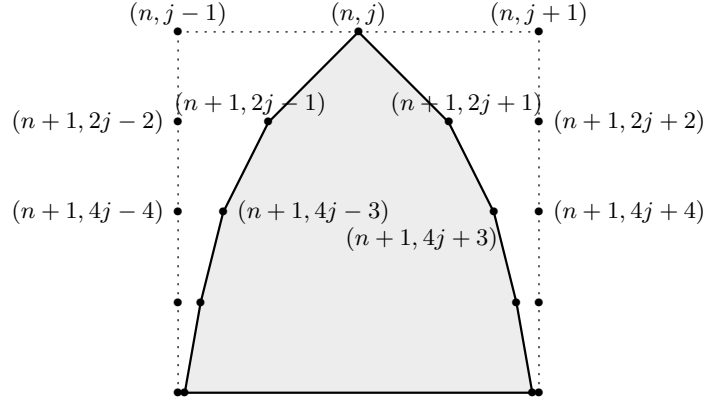

\begin{proof}
The inclusion $\II (\overline{E})\subseteq \II(E)$ is obvious by
\eqref{3.1}. We prove $\II(E)\subseteq I_x$ for all $x\in
\overline{E}$. Suppose ad absurdum there is $x\in \overline{E}$
for which $\II(E)\nsubseteq I_x$, i.e.~ there is $(n,j)\in\VV$
with $(n,j) \in D(\II(E))$ and $(n,j) \notin D(I_x)$. The latter
and Remark \ref{R3} yield
\begin{equation}\label{3.2}
r(n,j-1) <x<r(n,j+1) .
\end{equation}

On the other hand, because $D(\II(E))$ contains $(n,j)$, every
diagram $D(I_\theta)$, $\theta \in E$, must contain the whole
``pyramid" starting at $(n,j)$, see Figure \ref{Figure12}. Thus
\begin{equation*}
\forall \theta \in E,\ \forall k\geq 1, \quad \theta \in
[0,r(n+k,2^k j-2^k+1),1] \cup [r(n+k,2^k j+2^k-1),1].
\end{equation*}
But
\begin{equation*}
r(n+k,2^k j+2^k-1)=\frac{kp(n,j+1)+p(n,j)}{kq(n,j+1)+q(n,j)}\
\stackrel{k}{\longrightarrow} \ \frac{p(n,j+1)}{q(n,j+1)}
=r(n,j+1)
\end{equation*}
and
\begin{equation*}
r(n+k,2^k j-2^k+1)=\frac{kp(n,j-1)+p(n,j)}{kq(n,j-1)+q(n,j)}\
\stackrel{k}{\longrightarrow} \ \frac{p(n,j-1)}{q(n,j-1)}
=r(n,j-1),
\end{equation*}
hence
\begin{equation*}
E\subseteq [0,r(n,j-1)]\cup [r(n,j+1),1],
\end{equation*}
which is in contradiction with \eqref{3.2}.
\end{proof}

\begin{remark}\label{R8}
{\em We have $q(n,2j)=q(n-1,j) <\min \{ q(n,2j-1),q(n,2j+1) \}$,
so if $r(n,2j)=\frac{p}{q}$, then
\begin{equation*}
r(n,2j+1)-r(n,2j-1) =\frac{1}{q(n,2j-1) q(n,2j)}+\frac{1}{q(n,2j)
q(n,2j+1)} <\frac{2}{q^2} .
\end{equation*}
One can give a better estimate as follows. Let $\theta=\frac{p}{q}
\in (0,1)$ be a rational number in lowest terms and let
$\bar{p}\in \{ 1,\ldots,q-1\}$ denote the multiplicative inverse
of $p$ modulo $q$. Let $n_0=n_0(\theta)$ be the smallest $n$ such
that $\theta=r(n,j_0)$ for some $j_0$. Then $j_0$ is odd and the
labels $r^\prime=\frac{p^\prime}{q^\prime}$ and respectively
$r^{\prime\prime}=\frac{p^{\prime\prime}}{q^{\prime\prime}}$ of
the ``left parent" $(n_0-1,\frac{j_0-1}{2})$ and respectively of
the ``right parent" $(n_0-1,\frac{j_0+1}{2})$ of the vertex
$(n_0,j_0)$, are given by
$(q^\prime,p^\prime)=\big(\bar{p},\frac{p\bar{p}-1}{q}\big)$, and
respectively by
$(q^{\prime\prime},p^{\prime\prime})=\big(q-\bar{p},p-\frac{p\bar{p}-1}{q}\big)=(q,p)-(q^\prime,p^\prime)$.
Furthermore, we have $r(n_0+k,2^k
j_0-1)=\frac{p+kp^\prime}{q+kq^\prime}$, $r(n_0+k,2^k
j_0+1)=\frac{p+kp^{\prime\prime}}{q+kq^{\prime\prime}}$, and}
\begin{equation*}
\max \bigg\{ r(n_0+k,2^kj_0+1)-\frac{p}{q}\ ,\ \frac{p}{q}-
r(n_0+k,2^kj_0-1)\bigg\} <\frac{1}{kq^2}  .
\end{equation*}
\end{remark}

\begin{lemma}\label{L11}
For some $x\in [0,1]$ and $S\subseteq [0,1]$ suppose
$\II(S)\subseteq \II_x$. Then $x\in \overline{S}$.
\end{lemma}

\begin{proof}
Obviously two cases may occur:

{\em Case I: $x \notin \Q$.} Let $\big( \frac{p_n}{q_n}\big)$
denote the sequence of continued fraction approximations of $x$.
Taking stock on the definition of the ideal $\II_x$ we get
positive integers $k_1<k_2<\cdots$ and vertices $(k_n,j_n)\in
D(\AA)$ with the following properties:
\begin{itemize}
\item[(i)] $\quad \displaystyle r(k_n,j_{n})=\frac{p_n}{q_n}$;
\item[(ii)] $\quad j_{n}$ is even; \item[(iii)] $\quad (k_n,j_{n})
\notin D(\II_x)$.
\end{itemize}
Actually (iii) is a plain consequence of (i) and gives in turn,
cf.~ Remark \ref{R3},
\begin{equation}\label{3.3}
r(k_n,j_{n}-1)<x<r(k_n,j_{n}+1)  .
\end{equation}

{\em Case II: $x\in \Q$.} There is $n_0$ such that $(n,j_n) \notin
D(\II_x)$ and $r(n,j_n)=x$ for all $n\geq n_0$. In this case we
take $k_n=n$.

Suppose that $\exists n\geq n_0$, $\forall \theta \in S$,
$(k_n,j_{n})\in D(\II_{\theta})$. Then $(k_n,j_{n})\in
D(\II(S))\setminus D(\II_x)$, which contradicts the assumption of
the lemma. Therefore we must have
\begin{equation*}
\forall n,\ \exists \theta_n \in S, \ (k_n,j_{n}) \notin
D(\II_{\theta_n}),
\end{equation*}
which according to Remark \ref{R3} gives
\begin{equation}\label{3.4}
r(k_n,j_{n}-1) <\theta_n <r(k_n,j_{n}+1)  .
\end{equation}
From \eqref{3.3}, \eqref{3.4} and Remark \ref{R8} we now infer
\begin{equation*}
\vert x-\theta_n \vert <r(k_n,j_{n}+1)-r(k_n,j_{n}-1)
<\frac{2}{q_n^2} ,\qquad \forall n\geq n_0,
\end{equation*}
and so $\mathrm{dist} (x,S)=0$. This concludes the proof of the
lemma.
\end{proof}

As a consequence, the Jacobson topology is Hausdorff when
restricted to the subset $\Prim_0 \AA=\{ I_\theta :\theta \in
[0,1]\}$ of $\Prim \AA$. Moreover, we have

\begin{cor}\label{C12}
Let $( \theta_n)$ be a sequence in $[0,1]$. The following are
equivalent:
\begin{itemize}
\item[(i)] $\theta_n \rightarrow \theta$ in $[0,1]$. \item[(ii)]
$I_{\theta_n}\rightarrow I_\theta$ in $\Prim \AA$.
\end{itemize}
\end{cor}

\begin{proof}
(i) Suppose $\theta_n \rightarrow \theta$ in $[0,1]$ but
$I_{\theta_n} \nrightarrow I_\theta$ in $\Prim \AA$. Then there is
$I$ ideal in $\AA$ such that $I\nsubseteq I_\theta$ and there is a
subsequence $( n_k)$ such that $I_{\theta_{n_k}} \notin \OO_I$, so
that $I\subseteq I_{\theta_{n_k}}$. By Lemma \ref{L7} this also
yields $I\subseteq I_\theta$, which is a contradiction.

(ii) Suppose $I_{\theta_n}\rightarrow I_\theta$ in $\Prim \AA$ but
$\theta_n \nrightarrow \theta$ in $[0,1]$. Then there is a
subsequence $(n_k)$ such that $\theta \notin \overline{\{
\theta_{n_k}\}}$. By Lemma \ref{L11} we have $I:=\cap_k
I_{\theta_{n_k}} \nsubseteq I_\theta$, and so $I_\theta \in
\OO_I$. But on the other hand $I\subseteq I_{\theta_{n_k}}$, i.e.
$I_{\theta_{n_k}} \notin \OO_I$ for all $k$, thus contradicting
$I_{\theta_{n_k}} \rightarrow I_\theta$.
\end{proof}

\section{A description of the dimension group}
By a classical result of Elliott (\cite{Ell}, see also
\cite{Eff}), AF algebras are classified up to isomorphism by their
dimension groups. In this section we give a description of the
dimension group $K_0(\CCC)$ of the codimension one ideal
$\CCC=I_1$ of $\AA$ obtained by erasing all vertices $(n,2^n)$
from the Bratteli diagram. This is inspired by the generating
function identity \cite{Ca}
\begin{equation*}
\sum\limits_{n\geq 0} \theta_n X^n =\prod\limits_{k\geq 0}
(1+X^{2^k}+X^{2^{k+1}}),
\end{equation*}
where $( \theta_n)$ is the Stern-Brocot sequence
$q(0,0)$, $q(1,0)$, $q(1,1)$, $q(2,0)$, $q(2,1)$, $q(2,2)$,
$q(2,3)$, $\ldots$, $q(n,0)$, $\ldots$, $q(n,2^n-1)$, $q(n+1,0)$,
$\ldots$

For each integer $n\geq 0$, set
\begin{equation*}
p_{(n,k)}(X):=\begin{cases} 1 & \mbox{\rm if $k=0$,} \\ X^k+X^{-k}
& \mbox{\rm if $1\leq k<2^n$,}
\end{cases}
\end{equation*}
and consider the abelian additive group
\begin{equation*}
\PP_n:=\bigg\{ \sum\limits_{0\leq k<2^n} c_k p_{(n,k)}: c_k
\in\Z\bigg\}.
\end{equation*}
Set
\begin{equation*}
\varrho(X)=X^{-1}+1+X,\quad \varrho_n(X)=\prod\limits_{0\leq k<n}
\varrho (X^{2^k}),
\end{equation*}
and define the injective group morphisms
\begin{equation*}
\begin{split}
& \beta_m:\PP_m\rightarrow \PP_{m+1},\ \  (\beta_m (p))(X)=\varrho
(X)p(X^2),\\ & \beta_{m,n}:\PP_m \rightarrow \PP_n,\ \
(\beta_{m,n}(p))(X)=(\beta_{n-1}\cdots \beta_m
(p))(X)=\varrho_{m-n}(X)p(X^{2^{n-m}}),\quad m<n.
\end{split}
\end{equation*}
Note that
\begin{equation}\label{4.1}
\begin{split}
& (\beta_n (p_{(n,k)}))(X) =\varrho (X) p_{(n,k)}(X^2) \\ & \qquad
=\begin{cases} p_{(n+1,0)}(X)+p_{(n+1,1)}(X) & \mbox{\rm if $k=0$},\\
p_{(n+1,2k-1)}(X)+p_{(n+1,2k)}(X)+p_{(n+1,2k+1)}(X) & \mbox{\rm if
$1\leq k< 2^n$}.
\end{cases}
\end{split}
\end{equation}

The group $K_0(\CCC_n)$ identifies with the free abelian group
$\Z^{2^n}$, generated by the Murray-von Neumann equivalence
classes $[e_{(n,k)}]$ of minimal projections $e_{(n,k)}$ in the
central summand $\AA_{(n,k)}$, $0\leq k< 2^n$. We have
$K_0(\CCC)=\varinjlim (K_0(\CCC_n),\alpha_n)$, the injective
morphisms $\alpha_n :K_0(\CCC_n)\rightarrow K_0(\CCC_{n+1})$ being
given by
\begin{equation*}
\alpha_n ([e_{(n,k)}])=\begin{cases} [e_{(n+1,0)}]+[e_{(n+1,1)}] &
\mbox{\rm if $k=0$,}\\
[e_{(n+1,2k-1)}]+[e_{(n+1,2k)}]+[e_{(n+1,2k+1)}] & \mbox{\rm if
$1\leq k<2^n$.}
\end{cases}
\end{equation*}
The positive cone $K_0(\CCC_n)^+$ consists of elements of form
$\sum_{k=0}^{2^n-1} c_k [e_{(n,k)}]$, $c_k\in\Z_+$. The groups
$K_0(\CCC_n)$ and $\PP_n$ are identified by the group isomorphism
$\phi_n$ mapping $[e_{(n,k)}]$ onto $p_{(n,k)}$. Equalities
\eqref{4.1} are reflected into the commutativity of the diagram
\begin{equation}\label{4.2}
\xymatrix{ \mbox{$K_0(\CCC_n)$} \ar@{->} [r]^(.6){\phi_n}
\ar@{->} [d] _{\alpha_n} & \mbox{$\PP_n$} \ar@{->} [d] ^{\beta_n} \\
\mbox{$K_0(\CCC_{n+1})$} \ar@{->} [r] _(.6){\phi_{n+1}} &
\PP_{n+1} }
\end{equation}
As a result, $K_0(\CCC)$ is isomorphic with the abelian group
$\PP=\varinjlim (\PP_n,\beta_n)$ and can, therefore, be described
as $(\cup_n \PP_n)/_\sim =\Z [X+X^{-1}]/_\sim$ where $\sim$ is the
equivalence relation given by equality on each $\PP_n \times
\PP_n$, and for $p\in\PP_m$, $q\in\PP_n$, $m< n$, by
\begin{equation*}
p\sim q\ \Longleftrightarrow\ q(X)=(\beta_{m,n} (p) )
(X)=p(X^{2^{n-m}}) \prod\limits_{0\leq k<n-m}
(X^{-2^k}+1+X^{2^k}).
\end{equation*}
Let $[p]$ denote the equivalence class of $p\in \cup_n \PP_n$. The
addition on $\PP$ is given by
\begin{equation*}
[p]+[q]=[\beta_{m,n}(p)+q],\quad p\in\PP_m,\ q\in\PP_n,\ m\leq n,
\end{equation*}
and does not depend on the choice of $m$ or $n$. For example
\begin{equation*}
\begin{split}
[X^{-1}+X]+[X^{-3}+X^3] & =[(X^{-1}+1+X)(X^{-2}+X^2)+X^{-3}+X^3]
\\ & =[2(X^{-3}+X^3)+(X^{-2}+X^2)+(X^{-1}+X)].
\end{split}
\end{equation*}

An element $[p]$, $p\in\PP_n$, belongs to the positive cone
$\PP^+$ of the dimension group precisely when there is an integer
$N>n$ such that $\beta_{n,N} (p)$ has nonnegative coefficients.
The equality (where $c_{r+1}=0$)
\begin{equation*}
\begin{split}
& (X^{-1}+1+X)\sum\limits_{0\leq k<2^n} c_k (X^{2k}+X^{-2k})
\\ & \qquad \qquad =\sum\limits_{0\leq k<2^n} c_k (X^{2k}+X^{-2k}) +
\sum\limits_{0\leq k<2^n} (c_k+c_{k+1}) (X^{2k+1}+X^{-2k-1})
\end{split}
\end{equation*}
shows that $p(X)$ has nonnegative coefficients if and only if
$\varrho (X)p(X^2)$ has the same property. Therefore $[p]\in\PP^+$
precisely when $p(X)$ has nonnegative coefficients.

Consider the positive integers $q^\prime_{(n,k)}$, $n\geq 0$,
$0\leq k<2^n$, describing the sizes of central summands in
\begin{equation}\label{4.3}
\CCC_n =\bigoplus_{0\leq k<2^n} \M_{q^\prime_{(n,k)}},
\end{equation}
that is
\begin{equation*}
\begin{cases}
q^\prime_{(n,0)}=q^\prime_{(n,2^n-1)}=1 ,\\
q^\prime_{(n,2k)}=q^\prime_{(n-1,k)},\\
q^\prime_{(n,2k+1)}=q^\prime_{(n-1,k)}+q^\prime_{(n-1,k+1)},\quad
0\leq k<2^n .\end{cases}
\end{equation*}
For instance $q^\prime (3,k)$, $0\leq k\leq 7$, are given by
$1,3,2,3,1,2,1,1$, and $q^\prime (4,k)$, $0\leq k\leq 15$, by
$1,4,3,5,2,5,3,4,1,3,2,3,1,2,1,1$. From \eqref{4.3} we have
\begin{equation*}
\sum\limits_{0\leq k<2^n} q^\prime(n,k) [e_{(n,k)}] =[1] \quad
\mbox{\rm in $K_0(\CCC)$}.
\end{equation*}
This corresponds to
\begin{equation}\label{4.4}
\sum\limits_{0\leq k<2^n} q^\prime (n,k) p_{(n,k)}(X)=\varrho_n
(X).
\end{equation}

One can give a representation of $K_0(\CCC)$ where the injective
maps $\beta_n$ in \eqref{4.2} are replaced by inclusions $\iota_n
(p)=p$. Define
\begin{equation*}
\phi_{(n,k)}(X)
=\frac{p_{(n,k)}(X^{1/2^n})}{\varrho_{(n,k)}(X^{1/2^n})}
=\begin{cases} \displaystyle\frac{1}{\prod_{j=1}^n
(X^{-1/2^j}+1+X^{1/2^j})} & \mbox{\rm if $k=0$,}\\
\displaystyle\frac{X^{k/2^n}+X^{-k/2^n}}{\prod_{j=1}^n
(X^{-1/2^j}+1+X^{1/2^j})} & \mbox{\rm if $1\leq k<2^n,$} \\
\end{cases}
\end{equation*}
and consider the additive abelian group
\begin{equation*}
\RR_n :=\bigg\{ \sum\limits_{0\leq k<2^n} c_k \phi_{(n,k)}: c_k
\in \Z\bigg\} .
\end{equation*}
Equalities \eqref{4.1} become
\begin{equation*}
\begin{split}
\begin{cases} \phi_{(n+1,0)}+\phi_{(n+1,1)}=\phi_{(n,0)},\\
\phi_{(n+1,2k-1)}+\phi_{(n+1,2k)}+\phi_{(n+1,2k+1)}=\phi_{(n,k)},\quad
1\leq k<2^n ,
\end{cases}
\end{split}
\end{equation*}
and show that $\RR_n \subseteq \RR_{n+1}$ and that the diagram
$$
\xymatrix{ \mbox{$K_0(\CCC_n)$} \ar@{->} [r]^(.6){\psi_n}
\ar@{->} [d] _{\alpha_n} & \mbox{$\RR_n$} \ar@{^{(}->} [d] ^{\iota_n} \\
\mbox{$K_0(\CCC_{n+1})$} \ar@{->} [r] _(.6){\psi_{n+1}} &
\RR_{n+1} } $$ is commuting, where
$\psi([e_{(n,k)}])=\phi_{(n,k)}$. Therefore $K_0(\CCC)=\RR:=\cup_n
\RR_n$. Taking $X={\mathrm e}^Y$, we see that $K_0(\CCC)$ can be
viewed as the $\Z$-linear span of $\widetilde{\phi}_{(n,k)}$,
$n\geq 0$, $0\leq k<2^n$, where
\begin{equation*}
\widetilde{\phi}_{(n,k)}(Y)=\begin{cases} \vspace{0.1cm}
\displaystyle \frac{1}{\prod_{j=1}^n (1+2\cosh (Y/2^j))} &
\mbox{\rm if $k=0,$}\\ \displaystyle \frac{2\cosh
(kY/2^n)}{\prod_{j=1}^n (1+2\cosh (Y/2^j))}  & \mbox{\rm if $1\leq
k<2^n.$}
\end{cases} \qquad
\end{equation*}
One can certainly replace $Y$ by ${\mathrm i}Y$ and use $\cos$
instead of $\cosh$.

\section{Traces on $\AA$}
We augment the diagram $\GG=D(\AA)$ into $\widetilde{\GG}$, by
adding a $(-1)^{\mathrm{st}}$ floor with only one vertex
$\star=(-1,0)$ connected to both $(0,0)$ and $(0,1)$. Traces
$\tau$ on $\AA$ are in one-to-one correspondence (cf., e.g.,
\cite[Section 3.6]{EK}) with families
$\alpha^\tau=(\alpha^\tau_{(n,k)})$ of numbers in $[0,1]$, $n\geq
-1$, $0\leq k\leq 2^n$, such that
\begin{equation*}
\begin{cases}
\alpha^\tau_\star=1, & \\
\alpha^\tau_{(n,0)}=\alpha^\tau_{(n+1,0)}+\alpha^\tau_{(n+1,1)} &
\mbox{\rm if $n\geq -1 ,$}\\
\alpha^\tau_{(n,2^n)}=\alpha^\tau_{(n+1,2^{n+1})}+\alpha^\tau_{(n+1,2^{n+1}-1)}
& \mbox{\rm if $n\geq 0,$} \\
\alpha^\tau_{(n,k)}=\alpha^\tau_{(n+1,2k-1)}+\alpha^\tau_{(n+1,2k)}+\alpha^\tau_{(n+1,2k+1)}
& \mbox{\rm if $n\geq 1$, $0<k<2^n.$}
\end{cases}
\end{equation*}

An inspection of $\widetilde{\GG}$ shows that such a family
$\alpha^\tau$ is uniquely determined by the numbers
$\alpha^\tau_{(n,k)}$ with odd $k$. Let $\TT$ denote the diagram
obtained by removing the memory in $\widetilde{\GG}$. Its set of
vertices $V(\TT)$ consists of $\star$ and $(n,k)$ with $n\geq 0$
and odd $k$. For $v=(n,k)$ define $Lv=(n+1,2k-1)$ if $n\geq 0$,
$0< k\leq 2^n$, and $Rv=(n+1,2k+1)$ if $n\geq -1$, $0\leq k<2^n$.

\begin{figure}[htb]
\begin{center}
\unitlength 0.55mm
\begin{picture}(0,130)(0,-12)
%\thinlines
\path(-120,120)(120,100) \dottedline{2}(-120,120)(-120,100)
\dottedline{2}(-120,0)(-120,100)(0,80)
\dottedline{2}(120,0)(120,100)
\path(120,100)(0,80)(-60,60)(-90,40)(-105,20)(-112.5,0)
\dottedline{2}(0,80)(0,0)
\path(0,80)(60,60)(90,40)(105,20)(112.5,0)
\path(-60,60)(-30,40)(-15,20)(-7.5,0)
\path(60,60)(30,40)(15,20)(7.5,0) \dottedline{2}(-60,60)(-60,0)
\dottedline{2}(60,0)(60,60) \path(-90,40)(-75,20)(-67.5,0)
\path(-30,40)(-45,20)(-52.5,0) \path(90,40)(75,20)(67.5,0)
\path(30,40)(45,20)(52.5,0) \dottedline{2}(-90,40)(-90,0)
\dottedline{2}(-30,40)(-30,0) \dottedline{2}(90,40)(90,0)
\dottedline{2}(30,40)(30,0) \path(15,20)(22.5,0)
\path(45,20)(37.5,0) \path(75,20)(82.5,0) \path(105,20)(97.5,0)
\path(-15,20)(-22.5,0) \path(-45,20)(-37.5,0)
\path(-75,20)(-82.5,0) \path(-105,20)(-97.5,0)
\dottedline{2}(-120,80)(-60,60) \dottedline{2}(-120,60)(-90,40)
\dottedline{2}(-120,40)(-105,20) \dottedline{2}(-120,20)(-112.5,0)
\dottedline{2}(-97.5,0)(-90,20)(-82.5,0)\dottedline{2}(-67.5,0)(-60,20)(-52.5,0)
\dottedline{2}(-37.5,0)(-30,20)(-22.5,0)\dottedline{2}(-7.5,0)(0,20)(7.5,0)
\dottedline{2}(120,20)(112.5,0)
\dottedline{2}(97.5,0)(90,20)(82.5,0)\dottedline{2}(67.5,0)(60,20)(52.5,0)
\dottedline{2}(37.5,0)(30,20)(22.5,0)\dottedline{2}(120,80)(60,60)
\dottedline{2}(120,60)(90,40) \dottedline{2}(120,40)(105,20)
\dottedline{2}(105,20)(105,0)\dottedline{2}(75,20)(75,0)\dottedline{2}(45,20)(45,0)
\dottedline{2}(15,20)(15,0)
\dottedline{2}(-105,20)(-105,0)\dottedline{2}(-75,20)(-75,0)\dottedline{2}(-45,20)(-45,0)
\dottedline{2}(-15,20)(-15,0) \dottedline{2}(-30,40)(0,60)(30,40)
\dottedline{2}(-75,20)(-60,40)(-45,20)\dottedline{2}(75,20)(60,40)(45,20)
\dottedline{2}(-15,20)(0,40)(15,20)

\put(-125,100){\makebox(0,0){{\small $\mathbf{\frac{0}{1}}$}}}
\put(0,87){\makebox(0,0){{\small $\mathbf{\frac{1}{2}}$}}}
\put(125,100){\makebox(0,0){{\small $\mathbf{\frac{1}{1}}$}}}
\put(68,61){\makebox(0,0){{\small $\mathbf{\frac{2}{3}}$}}}
\put(-68,61){\makebox(0,0){{\small $\mathbf{\frac{1}{3}}$}}}

\put(-95,40){\makebox(0,0){{\small $\mathbf{\frac{1}{4}}$}}}
\put(-25,40){\makebox(0,0){{\small $\mathbf{\frac{2}{5}}$}}}
\put(25,40){\makebox(0,0){{\small $\mathbf{\frac{3}{5}}$}}}
\put(95,40){\makebox(0,0){{\small $\mathbf{\frac{3}{4}}$}}}
\put(-120,120){\makebox(0,0){{\large $\star$}}}

\put(-120,100){\makebox(0,0){{\tiny $\bullet$}}}
\put(120,100){\makebox(0,0){{\tiny $\bullet$}}}
\put(-120,80){\makebox(0,0){{\tiny $\bullet$}}}
\put(0,80){\makebox(0,0){{\tiny $\bullet$}}}
\put(120,80){\makebox(0,0){{\tiny $\bullet$}}}
\put(-120,60){\makebox(0,0){{\tiny $\bullet$}}}
\put(-60,60){\makebox(0,0){{\tiny $\bullet$}}}
\put(0,60){\makebox(0,0){{\tiny $\bullet$}}}
\put(60,60){\makebox(0,0){{\tiny $\bullet$}}}
\put(120,60){\makebox(0,0){{\tiny $\bullet$}}}

\put(-120,40){\makebox(0,0){{\tiny $\bullet$}}}
\put(-90,40){\makebox(0,0){{\tiny $\bullet$}}}
\put(-60,40){\makebox(0,0){{\tiny $\bullet$}}}
\put(-30,40){\makebox(0,0){{\tiny $\bullet$}}}
\put(0,40){\makebox(0,0){{\tiny $\bullet$}}}
\put(30,40){\makebox(0,0){{\tiny $\bullet$}}}
\put(60,40){\makebox(0,0){{\tiny $\bullet$}}}
\put(90,40){\makebox(0,0){{\tiny $\bullet$}}}
\put(120,40){\makebox(0,0){{\tiny $\bullet$}}}

\put(-110,20){\makebox(0,0){{\small $\mathbf{\frac{1}{5}}$}}}
\put(-70,20){\makebox(0,0){{\small $\mathbf{\frac{2}{7}}$}}}
\put(-50,20){\makebox(0,0){{\small $\mathbf{\frac{3}{8}}$}}}
\put(-10,20){\makebox(0,0){{\small $\mathbf{\frac{3}{7}}$}}}
\put(110,20){\makebox(0,0){{\small $\mathbf{\frac{4}{5}}$}}}
\put(70,20){\makebox(0,0){{\small $\mathbf{\frac{5}{7}}$}}}
\put(50,20){\makebox(0,0){{\small $\mathbf{\frac{5}{8}}$}}}
\put(10,20){\makebox(0,0){{\small $\mathbf{\frac{4}{7}}$}}}

\put(-120,20){\makebox(0,0){{\tiny $\bullet$}}}
\put(-105,20){\makebox(0,0){{\tiny $\bullet$}}}
\put(-90,20){\makebox(0,0){{\tiny $\bullet$}}}
\put(-75,20){\makebox(0,0){{\tiny $\bullet$}}}
\put(-60,20){\makebox(0,0){{\tiny $\bullet$}}}
\put(-45,20){\makebox(0,0){{\tiny $\bullet$}}}
\put(-30,20){\makebox(0,0){{\tiny $\bullet$}}}
\put(-15,20){\makebox(0,0){{\tiny $\bullet$}}}
\put(0,20){\makebox(0,0){{\tiny $\bullet$}}}
\put(120,20){\makebox(0,0){{\tiny $\bullet$}}}
\put(105,20){\makebox(0,0){{\tiny $\bullet$}}}
\put(90,20){\makebox(0,0){{\tiny $\bullet$}}}
\put(75,20){\makebox(0,0){{\tiny $\bullet$}}}
\put(60,20){\makebox(0,0){{\tiny $\bullet$}}}
\put(45,20){\makebox(0,0){{\tiny $\bullet$}}}
\put(30,20){\makebox(0,0){{\tiny $\bullet$}}}
\put(15,20){\makebox(0,0){{\tiny $\bullet$}}}

\put(-120,0){\makebox(0,0){{\tiny $\bullet$}}}
\put(-112.5,0){\makebox(0,0){{\tiny $\bullet$}}}
\put(-105,0){\makebox(0,0){{\tiny $\bullet$}}}
\put(-97.5,0){\makebox(0,0){{\tiny $\bullet$}}}
\put(-90,0){\makebox(0,0){{\tiny $\bullet$}}}
\put(-82.5,0){\makebox(0,0){{\tiny $\bullet$}}}
\put(-75,0){\makebox(0,0){{\tiny $\bullet$}}}
\put(-67.5,0){\makebox(0,0){{\tiny $\bullet$}}}
\put(-60,0){\makebox(0,0){{\tiny $\bullet$}}}
\put(-52.5,0){\makebox(0,0){{\tiny $\bullet$}}}
\put(-45,0){\makebox(0,0){{\tiny $\bullet$}}}
\put(-37.5,0){\makebox(0,0){{\tiny $\bullet$}}}
\put(-30,0){\makebox(0,0){{\tiny $\bullet$}}}
\put(-22.5,0){\makebox(0,0){{\tiny $\bullet$}}}
\put(-15,0){\makebox(0,0){{\tiny $\bullet$}}}
\put(-7.5,0){\makebox(0,0){{\tiny $\bullet$}}}
\put(0,0){\makebox(0,0){{\tiny $\bullet$}}}
\put(120,0){\makebox(0,0){{\tiny $\bullet$}}}
\put(112.5,0){\makebox(0,0){{\tiny $\bullet$}}}
\put(105,0){\makebox(0,0){{\tiny $\bullet$}}}
\put(97.5,0){\makebox(0,0){{\tiny $\bullet$}}}
\put(90,0){\makebox(0,0){{\tiny $\bullet$}}}
\put(82.5,0){\makebox(0,0){{\tiny $\bullet$}}}
\put(75,0){\makebox(0,0){{\tiny $\bullet$}}}
\put(67.5,0){\makebox(0,0){{\tiny $\bullet$}}}
\put(60,0){\makebox(0,0){{\tiny $\bullet$}}}
\put(52.5,0){\makebox(0,0){{\tiny $\bullet$}}}
\put(45,0){\makebox(0,0){{\tiny $\bullet$}}}
\put(37.5,0){\makebox(0,0){{\tiny $\bullet$}}}
\put(30,0){\makebox(0,0){{\tiny $\bullet$}}}
\put(22.5,0){\makebox(0,0){{\tiny $\bullet$}}}
\put(15,0){\makebox(0,0){{\tiny $\bullet$}}}
\put(7.5,0){\makebox(0,0){{\tiny $\bullet$}}}

\put(-112.5,-8){\makebox(0,0){{\small $\mathbf{\frac{1}{6}}$}}}
\put(-97.5,-8){\makebox(0,0){{\small $\mathbf{\frac{2}{9}}$}}}
\put(-82.5,-8){\makebox(0,0){{\small $\mathbf{\frac{3}{11}}$}}}
\put(-67.5,-8){\makebox(0,0){{\small $\mathbf{\frac{3}{10}}$}}}
\put(-52.5,-8){\makebox(0,0){{\small $\mathbf{\frac{4}{11}}$}}}
\put(-37.5,-8){\makebox(0,0){{\small $\mathbf{\frac{5}{13}}$}}}
\put(-22.5,-8){\makebox(0,0){{\small $\mathbf{\frac{5}{12}}$}}}
\put(-7.5,-8){\makebox(0,0){{\small $\mathbf{\frac{4}{9}}$}}}
\put(7.5,-8){\makebox(0,0){{\small $\mathbf{\frac{5}{9}}$}}}
\put(112.5,-8){\makebox(0,0){{\small $\mathbf{\frac{5}{6}}$}}}
\put(97.5,-8){\makebox(0,0){{\small $\mathbf{\frac{7}{9}}$}}}
\put(82.5,-8){\makebox(0,0){{\small $\mathbf{\frac{8}{11}}$}}}
\put(67.5,-8){\makebox(0,0){{\small $\mathbf{\frac{7}{10}}$}}}
\put(52.5,-8){\makebox(0,0){{\small $\mathbf{\frac{7}{11}}$}}}
\put(37.5,-8){\makebox(0,0){{\small $\mathbf{\frac{8}{13}}$}}}
\put(22.5,-8){\makebox(0,0){{\small $\mathbf{\frac{7}{12}}$}}}

\end{picture}
\end{center}
\caption{The diagram $\TT$} \label{Figure101}
\end{figure}
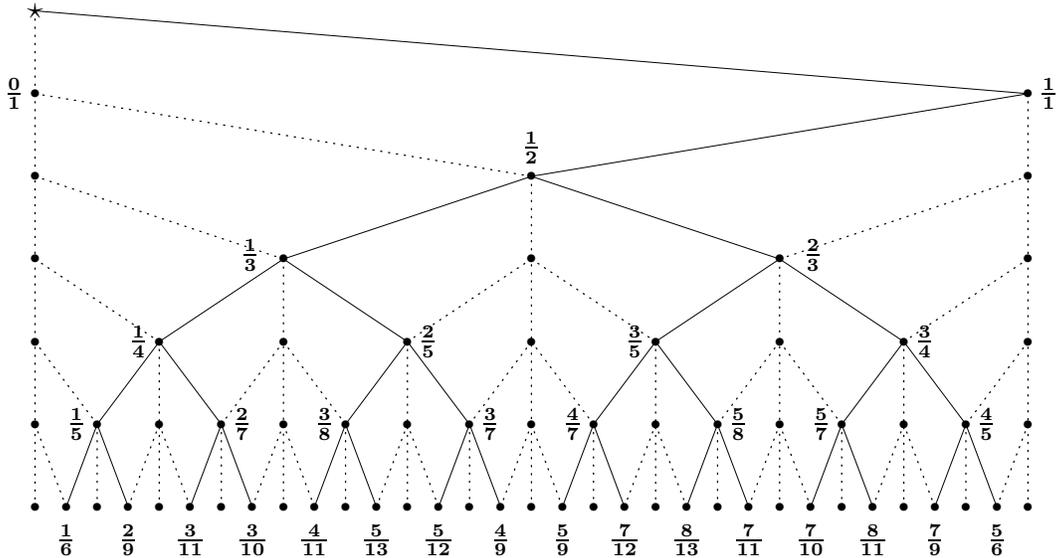

Given $\alpha^\tau_v$, $v=(n,k)\in V(\TT)$, define recursively for
$r\geq 1$
\begin{equation*}
\begin{cases}
\alpha^\tau_{(n+r,0)}=\alpha^\tau_{(n+r-1,0)}-\alpha^\tau_{(n+r,1)} &
\mbox{\rm if $n\geq -1$}, \\
\alpha^\tau_{(n+r,2^{n+r})}
=\alpha^\tau_{(n+r-1,2^{n+r-1})}-\alpha^\tau_{(n+r,2^{n+r}
-1)} & \mbox{\rm if $n\geq 0$},\\
\alpha^\tau_{(n+r,2^r k)} =\alpha^\tau_{(n+r-1,2^{r-1}k)}-
\alpha^\tau_{(n+r,2^r k-1)}-\alpha^\tau_{(n+r,2^r k+1)} &
\mbox{\rm if $n\geq 1$},
\end{cases}
\end{equation*}
or equivalently
\begin{equation}\label{5.1}
\begin{cases}
\alpha^\tau_{(n,0)}=\alpha^\tau_\star-\sum\limits_{j=0}^n \alpha^\tau_{(j,1)}
=\alpha^\tau_\star-\sum\limits_{j=0}^n \alpha^\tau_{L^j R\star} & \mbox{\rm if $n\geq 0$},\\
\alpha^\tau_{(n,2^n)} =\alpha^\tau_{(0,1)} -\sum\limits_{j=1}^n
\alpha^\tau_{(j,2^j-1)}=\alpha^\tau_{(0,1)}-\sum\limits_{j=1}^n
\alpha^\tau_{R^{j-1}L(0,1)}
 & \mbox{\rm if $n\geq 1$},\\
\alpha^\tau_{(n+r,2^r k)}=\alpha^\tau_{(n,k)}-\sum\limits_{j=1}^r
\big(
\alpha^\tau_{(n+j,2^{j}k-1)}+\alpha^\tau_{(n+j,2^{j}k+1)}\big) &
\\ \qquad \qquad \quad =\alpha^\tau_{(n,k)}-\sum\limits_{j=1}^r \big(
\alpha^\tau_{R^{j-1}L(n,k)} +\alpha^\tau_{L^{j-1}R(n,k)}\big) &
\mbox{\rm if $n\geq 2$}.
\end{cases}
\end{equation}

There is an obvious order relation on $V(\TT)$ defined by
$(n,k_n)\preceq (n^\prime,k_n^\prime)$ if $n\leq n^\prime$ and
there is a chain of vertices
$(n,k_n),\ldots,(n^\prime,k_n^\prime)$ such that $(n+i,k_{n+i})$
is connected to $(n+i+1,k_{n+i+1})$, i.e. $k_{n+i+1}-2k_{n+i}=\pm
1$. A function $f:V(\TT)\rightarrow \R$ is monotonically
decreasing if $f(v_1)\geq f(v_2)$ whenever $v_1 \preceq v_2$ in
$V(\TT)$. For each vertex $v=(n,k)\in V(\TT)$, let
\begin{equation}\label{5.2}
{\mathcal C}_v=\begin{cases} \{ L^j R\star:j\geq 0\} & \mbox{\rm
if $v=\star$,} \\  \{ R^{j-1}L(0,1):j\geq 1\} & \mbox{\rm if
$v=(0,1)$,} \\ \{ R^{j-1}Lv:j\geq 1\}\cup\{ L^{j-1}Rv:j\geq 1\} &
\mbox{\rm if $v\in V(\TT)\setminus \{ \star,(0,1)\}$,}
\end{cases}
\end{equation}
denote the set of vertices in $V(\TT)$ neighboring the vertical
infinite segment originating at $v$. As a result of \eqref{5.1}
and of non-negativity of $\alpha^\tau$ we have

\begin{proposition}\label{P5.1}
There is a one-to-one correspondence between traces on $\AA$ and
functions $\phi:V(\TT)\rightarrow [0,1]$ such that $\phi
(\star)=1$ and
\begin{equation}\label{5.3}
\phi(v)\geq \sum\limits_{w\in{\mathcal C}_v} \phi(w),\quad \forall
v\in V(\TT).
\end{equation}
\end{proposition}

Note that a function satisfying \eqref{5.3} is necessarily
monotonically decreasing.

\begin{figure}[htb]
\begin{center}
\unitlength 0.47mm
\begin{picture}(0,135)(0,-20)
%\thinlines

\path(-138,-15)(-134.5,-5) \path(-130.5,-5)(-127,-15)
\path(-108,-15)(-104.5,-5) \path(-100.5,-5)(-97,-15)
\path(-78,-15)(-72,-5) \path(-68,-5)(-62,-15)
\path(-38,-15)(-34.5,-5) \path(-32.5,-5)(-27,-15)

\put(-140,-20){\makebox(0,0){\scriptsize [6]}}
\put(-125,-20){\makebox(0,0){\scriptsize [4,2]}}
\put(-132.5,0){\makebox(0,0){\scriptsize [5]}}

\put(-110,-20){\makebox(0,0){\scriptsize [3,1,2]}}
\put(-95,-20){\makebox(0,0){\scriptsize [3,3]}}
\put(-102.5,0){\makebox(0,0){\scriptsize [3,2]}}
\put(-118.5,20){\makebox(0,0){\scriptsize [4]}}

\path(-130.5,5)(-120.5,15) \path(-116.5,15)(-104.5,5)

\put(-80,-20){\makebox(0,0){\scriptsize [2,1,3]}}
\put(-60,-20){\makebox(0,0){\scriptsize [2,1,1,2]}}
\put(-70,0){\makebox(0,0){\scriptsize [2,1,2]}}
\put(-40,-20){\makebox(0,0){\scriptsize [2,2,2]}}
\put(-25,-20){\makebox(0,0){\scriptsize [2,4]}}
\put(-32.5,0){\makebox(0,0){\scriptsize [2,3]}}
\put(-51.2,20){\makebox(0,0){\scriptsize [2,2]}}

\path(-68,5)(-53.2,15) \path(-49.2,15)(-34.5,5)
\path(-8,-15)(-2,-5) \path(2,-5)(8,-15)

\put(-10,-20){\makebox(0,0){\scriptsize [1,1,4]}}
\put(10,-20){\makebox(0,0){\scriptsize [1,1,2,2]}}
\put(0,0){\makebox(0,0){\scriptsize [1,1,3]}}

\path(37,-15)(45.5,-5) \path(49.5,-5)(58,-15)
\put(35,-20){\makebox(0,0){\scriptsize [1,1,1,1,2]}}
\put(60,-20){\makebox(0,0){\scriptsize [1,1,1,3]}}
\put(47.5,0){\makebox(0,0){\scriptsize [1,1,1,2]}}
\put(23.7,20){\makebox(0,0){\scriptsize [1,1,2]}}

\path(2,5)(21.7,15) \path(25.7,15)(45.5,5) \path(82,-15)(88,-5)
\path(92,-5)(98,-15) \put(80,-20){\makebox(0,0){\scriptsize
[1,2,3]}} \put(100,-20){\makebox(0,0){\scriptsize [1,2,1,2]}}
\put(90,0){\makebox(0,0){\scriptsize [1,2,2]}}

\path(122,-15)(125.5,-5) \path(129.5,-5)(133,-15)
\put(120,-20){\makebox(0,0){\scriptsize [1,3,2]}}
\put(135,-20){\makebox(0,0){\scriptsize [1,5]}}
\put(127.5,0){\makebox(0,0){\scriptsize [1,4]}}
\put(103.7,20){\makebox(0,0){\scriptsize [1,3]}}

\path(92,5)(101.7,15) \path(105.7,15)(125.5,5)
\put(-84.8,40){\makebox(0,0){\scriptsize [3]}}
\path(-116.5,25)(-86.8,35) \path(-82.8,35)(-53.2,25)
\put(63.7,40){\makebox(0,0){\scriptsize [1,2]}}
\path(25.7,25)(61.7,35) \path(65.7,35)(101.7,25)

\put(-10.5,60){\makebox(0,0){\scriptsize [2]}}
\path(-82.8,45)(-12.5,55) \path(-8.5,55)(61.7,45)
\put(145,80){\makebox(0,0){\scriptsize [1]}}
\path(-8.5,65)(143,75)

\put(-140,100){\makebox(0,0){\scriptsize 0}}
\path(-138,95)(143,85)
\end{picture}
\end{center}
\caption{The diagram $\TT$ in the continued fraction
representation} \label{Figure102}
\end{figure}
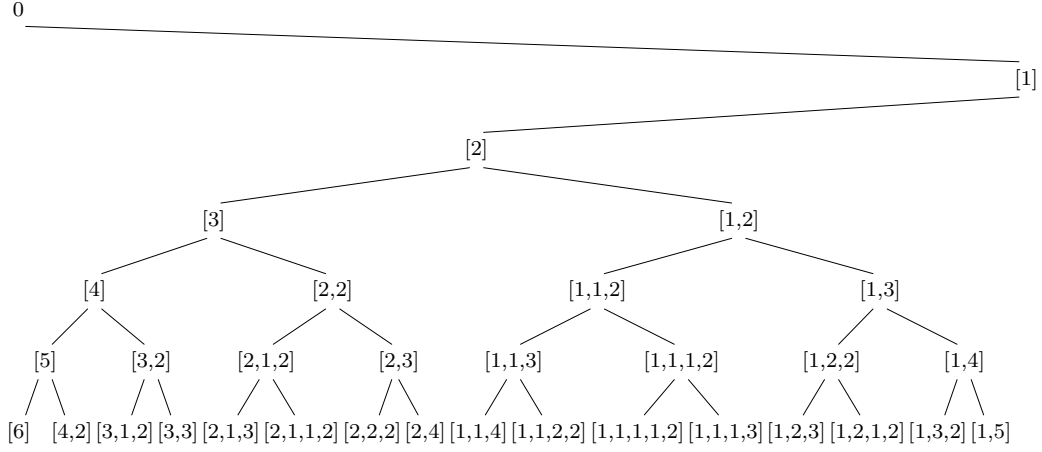

One can give a description of the set ${\mathcal C}_v$ using the
one-to-one correspondence $v\mapsto r(v)$ between the sets
$V(\TT)$ and $\Q\cap [0,1]$ (see Figure \ref{Figure102}). Any
number in $\Q\cap (0,1)$ can be uniquely represented as a
(reduced) continued fraction $[a_1,\ldots,a_t]$ with $a_t\geq 2$.
It is not hard to notice and prove that, for any $v\in V(\TT)$
with $r(v)=[a_1,\ldots,a_t]$, $a_t\geq 2$, we have
\begin{equation}\label{5.4}
\begin{split}
& r(Lv)=\begin{cases} [a_1,\ldots,a_{t-1},a_t -1,2] & \mbox{\rm if
$t$ even,}\\
[a_1,\ldots,a_{t-1},a_t+1] & \mbox{\rm if $t$ odd,}
\end{cases}
\\ & r(Rv)=\begin{cases} [a_1,\ldots,a_{t-1},a_t+1] & \mbox{\rm if $t$
even,} \\ [a_1,\ldots,a_{t-1},a_t -1,2] & \mbox{\rm if $t$ odd.}
\end{cases}
\end{split}
\end{equation}
As a result of \eqref{5.2} and \eqref{5.4} we have
\begin{equation*}
\{ r(w):w\in {\mathcal C}_v\}=\{ [a_1,\ldots,a_{t-1},a_t -1,1,k]:
k\geq 1\} \cup \{ [a_1,\ldots,a_{t-1},a_t,k]:k\geq 1\},
\end{equation*}
which shows in conjunction with Proposition \ref{P5.1} that there
is a  one-to-one correspondence between traces on $\AA$ and maps
$\phi :\Q \cap [0,1]\rightarrow [0,1]$ which satisfy
\begin{equation*}
\begin{cases}
1=\phi (0)\geq \sum\limits_{k=1}^\infty \phi
\big(\frac{1}{k}\big),\quad \phi (1)\geq \sum\limits_{k=1}^\infty
\phi\big(\frac{k}{k+1}\big),\\
\phi ([a_1,\ldots,a_t])\geq \sum\limits_{k=1}^\infty \big( \phi
([a_1,\ldots,a_{t-1},a_t -1,1,k])+\phi ([a_1,\ldots,a_{t-1},a_t,
k] \big),\quad a_t \geq 2.
\end{cases}
\end{equation*}

\section{Generators, relations, and braiding}
We shall use the path algebra model for AF algebras as in
\cite[Section 2.3.11]{GHJ} and \cite[Section 2.9]{EK}. Here
however a monotone increasing path $\xi$ will be encoded by the
sequence $(\xi_n)$ where $\xi_n$ gives the ``horizontal
coordinate" of the vertex at floor $n$, instead of its edges. To
use this model we again augment the diagram $\GG=D(\AA)$ into
$\widetilde{\GG}$.

Denote by $\Omega$ the (uncountable) set of monotone increasing
paths starting at $\star$. Let $\Omega_{[r}$ denote the set of
infinite monotone increasing paths starting on the
$r^{\mathrm{th}}$ floor of $\widetilde{\GG}$, $\Omega_{r]}$ the
set of monotone increasing paths that connect $\star$ with a
vertex on the $r^{\mathrm{th}}$ floor, and $\Omega_{[r,s]}$ the
set of monotone increasing paths starting on the $r^{\mathrm{th}}$
floor and ending on the $s^{\mathrm{th}}$ floor. Let
$\xi_r]\in\Omega_{r]}$, $\xi_{[r,s]}\in\Omega_{[r,s]}$,
$\xi_{[s}\in\Omega_{[s}$ denote the natural truncations of a path
$\xi \in \Omega$. By $\xi\circ\eta$ we denote the natural
concatenation of two paths $\xi\in\Omega_{r]}$ and
$\eta\in\Omega_{[r}$ with $\xi_{r}=\eta_{r}$. Consider the set
$R_r$ of pairs of paths $(\xi,\eta)\in \Omega_{r]}\times
\Omega_{r]}$ with the same endpoint $\xi_r=\eta_r$.  For each
$(\xi,\eta)\in R_r$ the mapping
\begin{equation*}
\Omega\ni\omega \mapsto T_{\xi,\eta} \omega
=\delta(\eta,\omega_{r]}) \xi\circ\omega_{[r}\in\Omega,
\end{equation*}
extends to a linear operator on the $\C$-linear space $\C \Omega$
with basis $\Omega$, and also to a bounded operator
$T_{\xi,\eta}:\ell^2(\Omega)\rightarrow \ell^2(\Omega)$ with $\|
T_{\xi,\eta}\|= 1$. We have $\AA=\overline{\cup_{r\geq 1}\AA_r}$
where the linear span $\AA_r$ of the operators $T_{\xi,\eta}$,
$(\xi,\eta)\in R_r$, forms a finite dimensional $C^*$-algebra as a
result of
\begin{equation*}
T_{\eta,\xi}^*=T_{\xi,\eta},\qquad
T_{\xi,\eta}T_{\xi^\prime,\eta^\prime}=\delta(\eta,\xi^\prime)
T_{\xi,\eta^\prime},\qquad \sum\limits_{\xi\in\Omega_{r]}}
T_{\xi,\xi}=1,
\end{equation*}
and the inclusion
$\AA_r\stackrel{\iota_r}{\hookrightarrow}\AA_{r+1}$ is given by
\begin{equation*}
\iota_r
(T_{\xi,\eta})=\sum\limits_{\substack{\lambda\in\Omega_{[r,r+1]}\\
\lambda_r=\xi_r (=\eta_r)}} T_{\xi\circ\lambda,\eta\circ\lambda}.
\end{equation*}

This model is employed to give a presentation by generators and
relations of the $C^*$-algebra $\AA$ in the spirit of the
presentation of the GICAR algebra from \cite[Example 2.23]{EK}. We
also construct two families of projections that satisfy
commutation relations reminiscent of the Temperley-Lieb relations.
Consider the following elements in $\AA$:
\begin{enumerate}
\item the projection $e_n$ in $\AA_{n-1,n} \subseteq \AA_n$ onto
the linear space of edges from N (north) to SW (south-west),
$n\geq 1$. \item the projection $f_n$ in $\AA_{n-1,n}\subseteq
\AA_n$ onto the linear span of edges from N to SE, $n\geq 0$.
\item the projection $g_n=1-e_n-f_n$ in $\AA_{n-1,n}\subseteq
\AA_n$ onto the linear span of edges from N to S, $n\geq 0$. \item
the partial isometry $v_n \in \AA_{n-1,n+1}\subseteq \AA_{n+1}$
with initial support $v_n^* v_n =\we_n =g_nf_{n+1}$ and final
support $v_nv_n^* =\wf_n =f_n e_{n+1}$, which flips paths in the
diamonds of shape N-S-SE-NE, $n\geq 0$. \item the partial isometry
$w_n \in \AA_{n-1,n+2}\subseteq \AA_{n+1}$ with initial support
$w_n^* w_n=\widetilde{e^\prime}_n =g_n e_{n+1}$ and final support
$w_nw_n^* =\wf_n^\prime =e_n f_{n+1}$, which flips paths in the
diamonds of shape N-S-SW-NW, $n\geq 1$.
\end{enumerate}
The AF-algebra $\AA$ is generated by the set ${\mathfrak G}=\{ e_n
\}_{n\geq 1} \cup \{ f_n \}_{n\geq 0} \cup \{ v_n \}_{n\geq 0}
\cup \{ w_n \}_{n\geq 1}$.

\begin{figure}[htb]
\begin{center}
\unitlength 0.45mm
\begin{picture}(0,155)(0,30)

\texture{c 0 c 0}
\shade\path(-120,120)(-120,150)(0,120)(-60,90)(-120,120)
\shade\path(-120,60)(-120,90)(-90,60)(-105,30)(-120,60)
\shade\path(-60,60)(-60,90)(-30,60)(-45,30)(-60,60)
\shade\path(0,60)(0,90)(30,60)(15,30)(0,60)
\shade\path(60,60)(60,90)(90,60)(75,30)(60,60)

\texture{c 000000}
\shade\path(0,120)(120,150)(120,120)(60,90)(0,120)
\shade\path(120,60)(120,90)(90,60)(105,30)(120,60)
\shade\path(60,60)(60,90)(30,60)(45,30)(60,60)
\shade\path(0,60)(0,90)(-30,60)(-15,30)(0,60)
\shade\path(-60,60)(-60,90)(-90,60)(-75,30)(-60,60)

\texture{ccc 00 ccc}
\shade\path(0,90)(0,120)(-60,90)(-30,60)(0,90)
\shade\path(60,90)(120,120)(120,90)(90,60)(60,90)
\shade\path(120,30)(120,60)(105,30)(120,30)
\shade\path(90,30)(90,60)(75,30)(90,30)
\shade\path(60,30)(60,60)(45,30)(60,30)
\shade\path(30,30)(30,60)(15,30)(30,30)
\shade\path(0,30)(0,60)(-15,30)(0,30)
\shade\path(-30,30)(-30,60)(-45,30)(-30,30)
\shade\path(-60,30)(-60,60)(-75,30)(-60,30)
\shade\path(-90,30)(-90,60)(-105,30)(-90,30)

\path(-120,150)(-120,180)(120,150)
\path(-120,120)(-120,150)(0,120)(120,150)(120,120)
\path(-120,120)(-120,90)
\path(-120,90)(-120,120)(-60,90)(0,120)(60,90)(120,120)(120,90)
\path(0,120)(0,90)

\path(-120,60)(-120,90)(-90,60)(-60,90)(-30,60)(0,90)(30,60)(60,90)(90,60)(120,90)(120,60)
\path(-60,90)(-60,60) \path(0,90)(0,60) \path(60,90)(60,60)
\path(-120,30)(-120,60)(-105,30)(-90,60)(-75,30)(-60,60)(-45,30)(-30,60)(-15,30)(0,60)(15,30)
(30,60)(45,30)(60,60)(75,30)(90,60)(105,30)(120,60)(120,30)
\path(-90,60)(-90,30) \path(-60,60)(-60,30) \path(-30,60)(-30,30)
\path(0,60)(0,30) \path(30,60)(30,30) \path(60,60)(60,30)
\path(90,60)(90,30)

\put(-120,180){\makebox(0,0){{\large $\star$}}}

\put(-5,145){\makebox(0,0){{$v_0$}}}
\put(-125,165){\makebox(0,0){{\small $g_0$}}}
\put(-125,135){\makebox(0,0){{\small $g_1$}}}
\put(125,135){\makebox(0,0){{\small $g_1$}}}
\put(-125,105){\makebox(0,0){{\small $g_2$}}}
\put(125,105){\makebox(0,0){{\small $g_2$}}}
\put(5,105){\makebox(0,0){{\small $g_2$}}}
\put(6,170){\makebox(0,0){{\small $f_0$}}}
\put(-55,140){\makebox(0,0){{\small $f_1$}}}
\put(55,140){\makebox(0,0){{\small $e_1$}}}
\put(-60,115){\makebox(0,0){{$v_1$}}}
\put(60,115){\makebox(0,0){{$w_1$}}}
\put(-85,110){\makebox(0,0){{\small $f_2$}}}
\put(35,110){\makebox(0,0){{\small $f_2$}}}
\put(85,110){\makebox(0,0){{\small $e_2$}}}
\put(-35,110){\makebox(0,0){{\small $e_2$}}}
\put(-90,90){\makebox(0,0){{$v_2$}}}
\put(-30,90){\makebox(0,0){{$w_2$}}}
\put(90,90){\makebox(0,0){{$w_2$}}}
\put(30,90){\makebox(0,0){{$v_2$}}}
\put(-105,60){\makebox(0,0){{$v_3$}}}
\put(-75,60){\makebox(0,0){{$w_3$}}}
\put(-45,60){\makebox(0,0){{$v_3$}}}
\put(-15,60){\makebox(0,0){{$w_3$}}}
\put(110,60){\makebox(0,0){{$w_3$}}}
\put(75,60){\makebox(0,0){{$v_3$}}}
\put(50,60){\makebox(0,0){{$w_3$}}}
\put(15,60){\makebox(0,0){{$v_3$}}}
\put(-125,70){\makebox(0,0){{\small $g_3$}}}
\put(-55,70){\makebox(0,0){{\small $g_3$}}}
\put(5,70){\makebox(0,0){{\small $g_3$}}}
\put(125,70){\makebox(0,0){{\small $g_3$}}}
\put(65,70){\makebox(0,0){{\small $g_3$}}}
\put(-97,75){\makebox(0,0){{\small $f_3$}}}
\put(-82,75){\makebox(0,0){{\small $e_3$}}}
\put(-37,75){\makebox(0,0){{\small $f_3$}}}
\put(-22,75){\makebox(0,0){{\small $e_3$}}}
\put(97,75){\makebox(0,0){{\small $e_3$}}}
\put(82,75){\makebox(0,0){{\small $f_3$}}}
\put(37,75){\makebox(0,0){{\small $e_3$}}}
\put(22,75){\makebox(0,0){{\small $f_3$}}}

\end{picture}
\end{center}
\caption{The generators of $\AA$} \label{Figure100}
\end{figure}
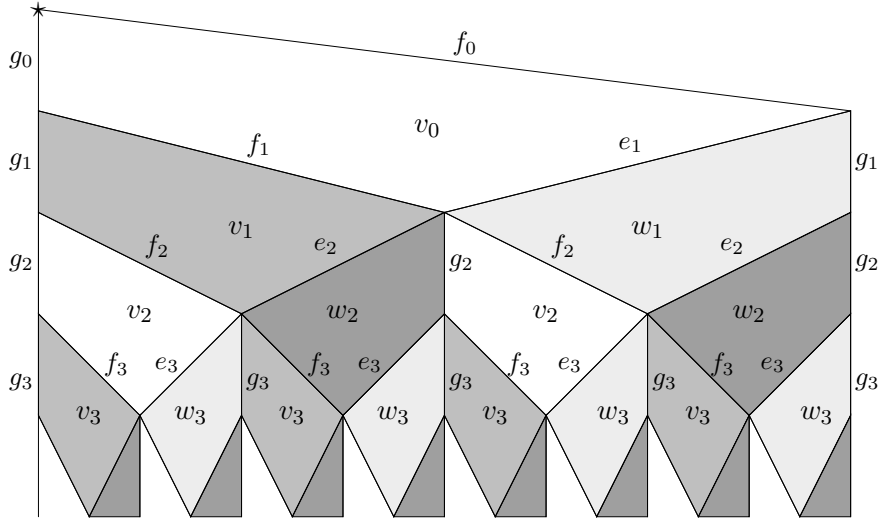

\begin{figure}[htb]
\begin{center}
\unitlength 0.4mm
\begin{picture}(70,40)(0,5)
%\thinlines
 \path(-80,10)(-50,40) \path(-20,10)(10,40) \path(40,10)(70,40)
\path(100,10)(130,40) \dottedline{3}(-50,10)(-50,40)(-20,10)
\dottedline{3}(10,10)(10,40)(40,10)
\dottedline{3}(70,10)(70,40)(100,10)
\dottedline{3}(130,10)(130,40)(160,10)
\dottedline{3}(-80,10)(-110,40)(-110,10)

\put(10,44){\makebox(0,0){{\footnotesize $(n-1,i)$}}}
\put(10,4){\makebox(0,0){{\footnotesize $(n,2i)$}}}
\put(-25,4){\makebox(0,0){{\footnotesize $(n,2i-1)$}}}
\put(45,4){\makebox(0,0){{\footnotesize $(n,2i+1)$}}}

\end{picture}
\end{center}
\caption{Support of projection $e_n$} \label{Figure13}
\end{figure}
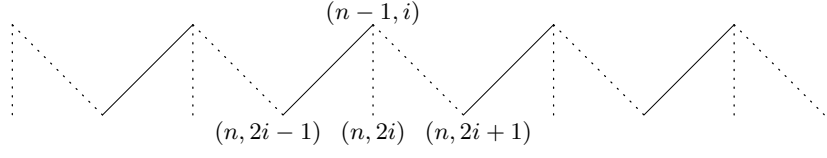

\begin{figure}[htb]
\begin{center}
\unitlength 0.4mm
\begin{picture}(70,40)(50,5)
\dottedline{3}(-50,40)(-50,10)
\dottedline{3}(-20,10)(10,40)(10,10)
\dottedline{3}(40,10)(70,40)(70,10)
\dottedline{3}(160,10)(190,40)(190,10)
\dottedline{3}(100,10)(130,40)(130,10) \path(-50,40)(-20,10)
\path(10,40)(40,10) \path(70,40)(100,10) \path(130,40)(160,10)

\put(70,44){\makebox(0,0){{\footnotesize $(n-1,i)$}}}
\put(70,4){\makebox(0,0){{\footnotesize $(n,2i)$}}}
\put(35,4){\makebox(0,0){{\footnotesize $(n,2i-1)$}}}
\put(105,4){\makebox(0,0){{\footnotesize $(n,2i+1)$}}}

\end{picture}
\end{center}
\caption{Support of projection $f_n$} \label{Figure14}
\end{figure}
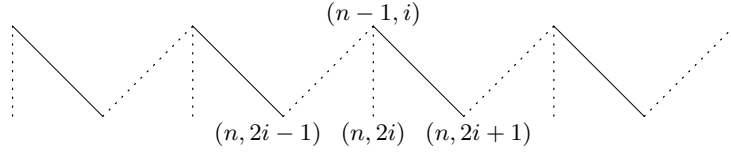

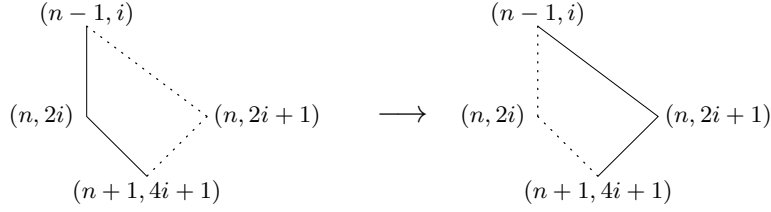
\begin{figure}[htb]
\begin{center}
\unitlength 0.4mm
\begin{picture}(70,55)(0,-10)
\dottedline{3}(-60,40)(-20,10)(-40,-10)
\path(-60,40)(-60,10)(-40,-10)

\path(90,40)(130,10)(110,-10)
\dottedline{3}(90,40)(90,10)(110,-10)

\put(-60,44){\makebox(0,0){{\footnotesize $(n-1,i)$}}}
\put(-75,10){\makebox(0,0){{\footnotesize $(n,2i)$}}}
\put(-40,-15){\makebox(0,0){{\footnotesize $(n+1,4i+1)$}}}
\put(0,10){\makebox(0,0){{\footnotesize $(n,2i+1)$}}}
\put(45,10){\makebox(0,0){{\large $\longrightarrow$}}}
\put(90,44){\makebox(0,0){{\footnotesize $(n-1,i)$}}}
\put(75,10){\makebox(0,0){{\footnotesize $(n,2i)$}}}
\put(110,-15){\makebox(0,0){{\footnotesize $(n+1,4i+1)$}}}
\put(150,10){\makebox(0,0){{\footnotesize $(n,2i+1)$}}}

\end{picture}
\end{center}
\caption{The partial isometry $v_n:g_n f_{n+1} \mapsto f_n
e_{n+1}$} \label{Figure15}
\end{figure}

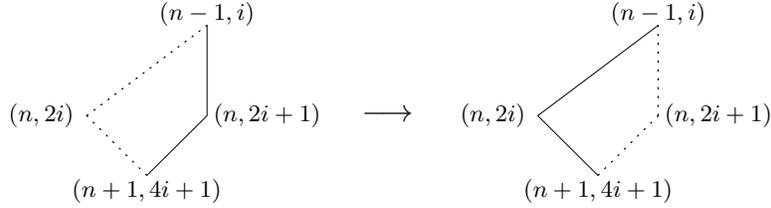
\begin{figure}[htb]
\begin{center}
\unitlength 0.4mm
\begin{picture}(70,55)(0,-10)
\dottedline{3}(-20,40)(-60,10)(-40,-10)
\path(-20,40)(-20,10)(-40,-10) \path(130,40)(90,10)(110,-10)
\dottedline{3}(130,40)(130,10)(110,-10)

\put(-20,44){\makebox(0,0){{\footnotesize $(n-1,i)$}}}
\put(-75,10){\makebox(0,0){{\footnotesize $(n,2i)$}}}
\put(-40,-15){\makebox(0,0){{\footnotesize $(n+1,4i+1)$}}}
\put(0,10){\makebox(0,0){{\footnotesize $(n,2i+1)$}}}
\put(40,10){\makebox(0,0){{\large $\longrightarrow$}}}
\put(130,44){\makebox(0,0){{\footnotesize $(n-1,i)$}}}
\put(75,10){\makebox(0,0){{\footnotesize $(n,2i)$}}}
\put(110,-15){\makebox(0,0){{\footnotesize $(n+1,4i+1)$}}}
\put(150,10){\makebox(0,0){{\footnotesize $(n,2i+1)$}}}

\end{picture}
\end{center}
\caption{The partial isometry $w_n:g_n e_{n+1}\mapsto e_n
f_{n+1}$} \label{Figure16}
\end{figure}

Straightforward commutation relations arise since elements defined
by edges that reach up to floor $\leq r$ commute with elements
defined by edges between the $r^{\mathrm{th}}$ and the
$s^{\mathrm{th}}$ floors with $r<s$, as a result of
$[\AA_r,\AA_r^\prime \cap \AA_s]=0$. For instance $v_s$ commutes
with $e_r,f_r,g_r$ if $r\leq s-1$ or $r\geq s+2$, and
$[v_s,v_r]=[v_s,v_r^*]=[v_s,w_r]=[v_s,w_r^*]=0$ if $\vert r-s\vert
\geq 2$. Besides, the elements of ${\mathfrak G}$ satisfy the
following commutation relations:
\begin{enumerate}
\item[(R1)] $e_n^2=e_n^*=e_n$, $\ f_n^2=f_n^*=f_n$, $\
g_n^2=g_n^*=g_n$, $\ e_n+f_n+g_n=1$;

$e_n,f_m,g_k$ mutually commute. \item[(R2)]
$(1-f_n)v_n=(1-e_{n+1})v_n=0$,
$v_n(1-g_n)=v_n(1-f_{n+1})=0$.\\
$(1-e_n)w_n=(1-f_{n+1})w_n=0$, $w_n(1-g_n)=w_n(1-e_{n+1})=0$.
\item[(R3)] $v_n g_n=f_nv_n$, $v_n f_{n+1}=e_{n+1}v_n$, $w_n
g_n=e_n w_n$, $w_n e_{n+1}=f_{n+1} w_n$. \item[(R4)] $v_n^*
v_n=g_n f_{n+1}$, $v_n v_n^*=f_n e_{n+1}$, $ w_n^* w_n=g_n
e_{n+1}$, $w_n w_n^*=e_n f_{n+1}$.
\end{enumerate}

As a result of (R1)--(R4) we also get
\begin{equation}\label{6.1}
\begin{aligned}
& v_{n+1}v_n=v_n^2= v_{n\pm 1} v_n^*=v_{n\pm 1}^* v_n=0,\\ &
w_{n+1}w_n=w_n^2 =w_{n\pm 1} w_n^*=w_{n\pm 1}^* w_n=0 ,\\
& v_nw_n=v_{n\pm 1} w_n=w_n v_n=w_{n\pm 1}v_n=0,\\
& v_nw_n^*=v_{n\pm 1} w_n^*=v_n^* w_n=v_n^* w_{n-1}=0.
\end{aligned}
\end{equation}
The only non-zero products $ab$ with $a\in \{ v_n,v_n^*
,w_n,w_n^*\}$ and $b\in\{ v_{n+1},v_{n+1}^*, w_{n+1},w_{n+1}^* \}$
are $v_n v_{n+1}$, $w_nw_{n+1}$, $w_n^*v_{n+1}$, and
$v_n^*w_{n+1}$.

\begin{figure}[htb]
\begin{center}
\unitlength 0.4mm
\begin{picture}(0,130)(0,15)
%\thinlines
\path(-80,140)(-80,100)(-60,80)\dottedline{3}(-80,140)(-20,120)(-60,80)

\path(-60,10)(-80,30)(-80,50)(-20,70)
\dottedline{3}(-20,70)(-20,50)(-40,30)(-60,10)

\put(-60,10){\makebox(0,0){{\tiny $\bullet$}}}
\put(-80,30){\makebox(0,0){{\tiny $\bullet$}}}
\put(-80,50){\makebox(0,0){{\tiny $\bullet$}}}
\put(-20,70){\makebox(0,0){{\tiny $\bullet$}}}
\put(-20,50){\makebox(0,0){{\tiny $\bullet$}}}
\put(-40,30){\makebox(0,0){{\tiny $\bullet$}}}

\dottedline{3}(120,140)(60,120)(80,100)(100,80)
\path(120,140)(120,120)(120,100)(100,80)

\put(120,140){\makebox(0,0){{\tiny $\bullet$}}}
\put(120,120){\makebox(0,0){{\tiny $\bullet$}}}
\put(120,100){\makebox(0,0){{\tiny $\bullet$}}}
\put(100,80){\makebox(0,0){{\tiny $\bullet$}}}
\put(80,100){\makebox(0,0){{\tiny $\bullet$}}}
\put(60,120){\makebox(0,0){{\tiny $\bullet$}}}

\put(-80,140){\makebox(0,0){{\tiny $\bullet$}}}
\put(-80,120){\makebox(0,0){{\tiny $\bullet$}}}
\put(-80,100){\makebox(0,0){{\tiny $\bullet$}}}
\put(-60,80){\makebox(0,0){{\tiny $\bullet$}}}
\put(-40,100){\makebox(0,0){{\tiny $\bullet$}}}
\put(-20,120){\makebox(0,0){{\tiny $\bullet$}}}

\put(-110,115){\makebox(0,0){$v_nv_{n+1}=$}}
\put(30,115){\makebox(0,0){$w_nw_{n+1}=$}}

\put(-110,40){\makebox(0,0){$w_n^* v_{n+1}=$}}
\put(35,40){\makebox(0,0){$v_{n}^* w_{n+1}=$}}

\dottedline{3}(60,70)(60,50)(80,30)(100,10)
\path(100,10)(120,30)(120,50)(60,70)

\put(60,70){\makebox(0,0){{\tiny $\bullet$}}}
\put(60,50){\makebox(0,0){{\tiny $\bullet$}}}
\put(80,30){\makebox(0,0){{\tiny $\bullet$}}}
\put(100,10){\makebox(0,0){{\tiny $\bullet$}}}
\put(120,30){\makebox(0,0){{\tiny $\bullet$}}}
\put(120,50){\makebox(0,0){{\tiny $\bullet$}}}
\end{picture}
\end{center}
\caption{The partial isometries $v_n v_{n+1}:g_n g_{n+1}
f_{n+2}\mapsto f_n e_{n+1} e_{n+2}$, $w_nw_{n+1}:g_n g_{n+1}
e_{n+2} \mapsto e_n f_{n+1} f_{n+2}$, $w_n^*v_{n+1}: e_n
g_{n+1}f_{n+2} \mapsto g_n e_{n+1} e_{n+2}$, $v_n^*w_{n+1}:f_n
g_{n+1}e_{n+2} \mapsto g_n f_{n+1}f_{n+2}$} \label{Figure18}
\end{figure}
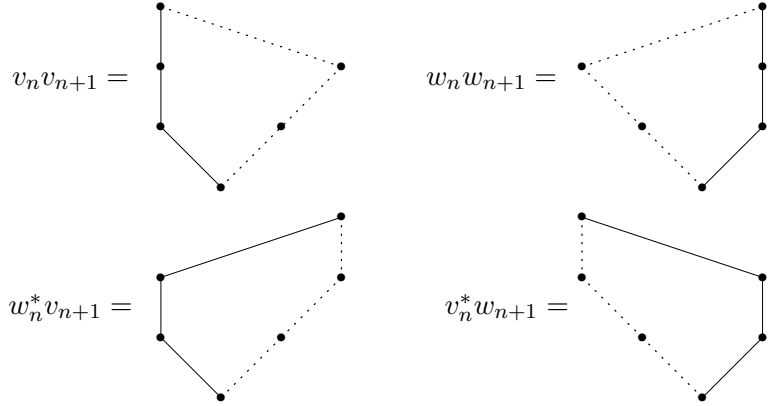

Let $B_n$ denote Artin's braid group generated by
$\sigma_1,\ldots,\sigma_{n-1}$ with relations
$\sigma_i\sigma_j=\sigma_j\sigma_i$ if $\vert i-j\vert >1$ and
$\sigma_i \sigma_{i+1} \sigma_i=\sigma_{i+1} \sigma_i
\sigma_{i+1}$. Relations \eqref{6.1} show in particular that
the partial isometries $v_{i-1}$, respectively $w_i$, satisfy
these braid relations.

Taking $R_n(\lambda):=1+\lambda v_n$, the equalities
\begin{equation}\label{6.3}
v_n^2=0,\qquad v_n v_{n\pm 1}v_n=0,
\end{equation}
yield the Yang-Baxter type relation
\begin{equation}\label{6.4}
R_n(\lambda)R_{n+1}(\lambda+\mu)R_n(\mu)=R_{n+1}(\mu)R_n(\lambda+\mu)R_{n+1}(\lambda).
\end{equation}

By analogy with the construction of Temperley-Lieb-Jones
projections in the GICAR algebra (cf., e.g., \cite{EK} or
\cite{GHJ}) for each $\lambda >0$ we put
$\tau=\frac{\lambda}{(1+\lambda)^2}\in ( 0,\frac{1}{4}]$ and
consider
\begin{equation}\label{6.5}
E_n =\frac{1}{1+\lambda}\ \big( v_n^*v_n+\sqrt{\lambda}\,
v_n+\sqrt{\lambda}\, v_n^*+\lambda v_nv_n^* \big)\in \AA,\quad
n\geq 0,
\end{equation}
\begin{equation}\label{6.6}
F_n =\frac{1}{1+\lambda} \ \big( w_n^*w_n+\sqrt{\lambda}\
w_n+\sqrt{\lambda}\ w_n^* +\lambda w_nw_n^*\big)\in\AA ,\quad
n\geq 1.
\end{equation}

\begin{proposition}
The elements $E_n$ and $F_n$ define {\em (}self-adjoint{\em )}
projections in the AF algebra $\AA$ satisfying the braiding
relations
\begin{equation}\label{6.7}
E_nF_n=F_nE_n=0,
\end{equation}
\begin{equation}\label{6.8}
[E_n,E_m]=[F_n,F_m]=[E_n,F_m]=0 \quad \mbox{if $\vert n-m\vert
\geq 2$,}
\end{equation}
\begin{equation}\label{6.9}
E_n E_{n+1}E_n =\tau E_n e_{n+2},\quad E_{n+1}E_nE_{n+1}=\tau
E_{n+1} g_n,
\end{equation}
\begin{equation}\label{6.10}
F_n F_{n+ 1} F_n =\tau F_n f_{n+2},\quad F_{n+1}F_nF_{n+1}=\tau
F_{n+1} g_n,
\end{equation}
\begin{equation}\label{6.11}
E_n F_{n+1}E_n =\lambda\tau E_n f_{n+2},\quad F_n E_{n+1}F_n
=\lambda\tau F_n e_{n+2},
\end{equation}
\begin{equation}\label{6.12}
E_{n+1}F_n E_{n+1}=\lambda\tau E_{n+1}e_n,\quad
F_{n+1}E_nF_{n+1}=\lambda\tau F_{n+1}f_n,
\end{equation}
\begin{equation}\label{6.13}
E_n E_{n+1}F_n=E_n F_{n+1}F_n=E_{n+1}E_n F_{n+1}=E_{n+1}F_n
F_{n+1}=0,
\end{equation}
\begin{equation}\label{6.14}
F_n E_{n+1}E_n=F_n F_{n+1} E_n =F_{n+1} E_n E_{n+1}=F_{n+1}
F_nE_{n+1} =0.
\end{equation}
\end{proposition}

\begin{proof}
The initial and final projections of the partial isometry $v_n$
are orthogonal, thus $E_n$ defines a projection in $\AA_n$ for
every $\lambda \geq 0$. A similar property holds for $F_n$, which
is seen to be orthogonal to $E_n$. The commutation relations
\eqref{6.8} are obvious because $v_{n+2}$ and $w_{n+2}$ commute
with all elements in $\AA_{n+1}$, including $E_n$ and $F_n$. By
\eqref{6.1} we have $v_n^*E_{n+1}=v_n v_{n+1}^*=0$, leading to
\begin{equation}\label{6.15}
E_n E_{n+1} =\frac{\sqrt{\lambda}}{(1+\lambda)^2}\ \big( v_n^* v_n
+\sqrt{\lambda}\ v_n\big)\big(v_{n+1}+\sqrt{\lambda}\
v_{n+1}v_{n+1}^*\big),
\end{equation}
and also
\begin{equation}\label{6.16}
E_{n+1}E_n=(E_nE_{n+1})^*=\frac{\sqrt{\lambda}}{(1+\lambda)^2}\
\big(v_{n+1}^*+\sqrt{\lambda}\ v_{n+1} v_{n+1}^*\big)\big( v_n^*
v_n +\sqrt{\lambda}\ v_n^*\big).
\end{equation}
From \eqref{6.15} and $v_{n+1}E_n=v_{n+1}^* v_n=0$ we have
\begin{equation}\label{6.17}
\begin{split}
E_nE_{n+1}E_n & =\frac{\lambda}{(1+\lambda)^3}\ \big(
v_n^*v_n+\sqrt{\lambda}\ v_n\big) v_{n+1}v_{n+1}^* \big( v_n^*
v_n+\sqrt{\lambda}\ v_n^*\big) \\ &
=\frac{\lambda}{(1+\lambda)^3}\ \big( \we_n+\sqrt{\lambda}\ v_n
\big) \wf_{n+1} \big( \we_n+\sqrt{\lambda}\ v_n^*\big).
\end{split}
\end{equation}
But $\we_n\wf_{n+1}\we_n=\we_n\wf_{n+1}=g_nf_{n+1}e_{n+1}=\we_n
e_{n+2}$, $v_n\wf_{n+1}\we_n=v_n \we_n e_{n+1}e_{n+2}=v_n e_{n+2}$
(and because $[e_{n+2},v_n]=0$ this also gives
$\we_n\wf_{n+2}v_n^*=v_n^*e_{n+2}$), and $v_n \wf_{n+1} v_n^*=v_n
f_{n+1}e_{n+2}v_n^* =v_n f_{n+1}v_n^* e_{n+2}=v_n g_n f_{n+1}v_n^*
e_{n+2}=v_n v_n^*e_{n+2}$, which we insert in \eqref{6.17} to get
\begin{equation*}
E_nE_{n+1}E_n =\tau E_n e_{n+2} .
\end{equation*}

From \eqref{6.16} and $v_n^*E_{n+1}=v_n^*v_{n+1}^*=0$ we find
\begin{equation}\label{6.18}
E_{n+1}E_nE_{n+1}=\frac{\lambda}{(1+\lambda)^3}\ \big( v^*_{n+1}
+\sqrt{\lambda}\ \wf_{n+1}\big) \we_n \big(
v_{n+1}+\sqrt{\lambda}\ \wf_{n+1}\big).
\end{equation}
As a result of $[g_n,v_{n+1}]=0$ and $(1-f_{n+1})v_{n+1}=0$ we
have $v^*_{n+1}\we_n v_{n+1}=\we_{n+1}g_n$. It is also plain that
$\wf_{n+1}\we_n\wf_{n+1}=\wf_{n+1}\we_n=\wf_{n+1}g_n$,
$\wf_{n+1}\we_nv_{n+1}=\wf_{n+1}g_nv_{n+1}=\wf_{n+1}v_{n+1}g_n=v_{n+1}g_n$,
and $v^*_{n+1}\we_n \wf_{n+1}=v^*_{n+1} \wf_{n+1} g_n=v^*_{n+1}
g_n$. Together with \eqref{6.18} these equalities yield
\begin{equation*}
E_{n+1}E_nE_{n+1} =\tau E_{n+1} g_n.
\end{equation*}

Equalities \eqref{6.10}--\eqref{6.13} are checked in a similar
way. \eqref{6.14} follows by taking adjoints in \eqref{6.13}.
\end{proof}

\section*{Acknowledgments}
I am grateful to Ola Bratteli, Marius Dadarlat, George Elliott,
Andreas Knauf, and Bruce Reznick for useful comments and
suggestions.

\end{document}